%% file: manuscriptR2.tex
\documentclass[a4paper, 11pt]{article}

\usepackage[margin=1.2in]{geometry}
\usepackage{parskip}
\usepackage[titletoc,title]{appendix}

\usepackage{array}
\usepackage{bm}
\usepackage{caption}
\usepackage{rotating}
\usepackage[hidelinks, colorlinks = true, linkcolor = blue, urlcolor = blue, anchorcolor = blue, citecolor = blue]{hyperref}
\usepackage{blkarray, bigstrut}
\usepackage{multirow}
\usepackage{booktabs}
\usepackage{siunitx}
\usepackage{mathtools}
\usepackage{amsmath}
\usepackage{amssymb}
\usepackage{amsthm}
\usepackage{complexity}
\usepackage{natbib}
\usepackage{graphicx}
\usepackage{pgf}
\usepackage{tikz}
\usetikzlibrary{matrix,patterns,positioning,decorations,arrows,shapes,decorations.markings,calc}
\usepackage[latin1]{inputenc}
\usepackage{verbatim}
\usepackage{subfiles}
\usepackage{pgfplots}
\usepackage{color}
\usepackage{float}
\usepackage{enumerate}
\usepackage{emptypage}
\usepackage{subcaption}
\usepackage{setspace}
\usepackage[export]{adjustbox}
\usepackage{subfiles}
\usepackage[algoruled]{algorithm2e}
\usepackage{xcolor}
\usepackage[textsize=normal]{todonotes}

\definecolor{borange50}{RGB}{254,243,230}
\definecolor{borange200}{RGB}{249,193,120}
\definecolor{borange600}{RGB}{186,100,14}
\definecolor{borange800}{RGB}{120,62,6}

\definecolor{bgray200}{RGB}{190,190,185}
\definecolor{bgray400}{RGB}{136,135,128}
\definecolor{bgray600}{RGB}{95,94,90}
\definecolor{bgray800}{RGB}{50,49,47}

\title{Branch and Price for Railway Crew Scheduling: Benchmark Instances and Computational Study}

\author{B.T.C. van Rossum
\vspace{0.1cm}\\
\small{Department of Industrial Engineering \& Innovation Sciences}\\ 
\small{Eindhoven University of Technology, the Netherlands}
\vspace{0.1cm}\\
\small{b.t.c.v.rossum@tue.nl}
}

\date{}

\begin{document}
  
\maketitle

\begin{abstract}
\noindent 
Railway crew scheduling consists of assigning a set of tasks to crew members in the form of feasible duties so as to minimise costs, and is one of the key planning problems faced by railway operators. While column generation is the dominant solution method, literature on exact branch-and-price approaches is scarce. Moreover, few publicly available benchmark instances exist. We conduct a systematic computational study of branch and price for railway crew scheduling, evaluating the contribution of branching rules, primal heuristics, and reduced cost fixing. To this end, we introduce a novel library of 100 benchmark instances, ranging in size from 90 to 3,016 tasks, each provided in a single-depot and a multi-depot variant. The instances are generated by mimicking the railway planning process, validated against data from Netherlands Railways, and publicly available with best known solutions and a solution checker. We report average optimality gaps below 1\% across all instance sizes, driven primarily by the strong performance of primal heuristics, whereas branching and reduced cost fixing contribute little to closing the gap. Branch and price also outperforms a multi-commodity flow formulation, which fails to find a feasible solution on many of the smallest instances. Our findings provide support for the popularity of column generation-based primal heuristics, and show that improving the dual bound remains the main challenge for exact methods in railway crew scheduling.
\vspace{5mm}
\end{abstract}

\section{Introduction}
\label{sec:intro}

Railway crew scheduling, the problem of assigning train trips to crew members, is one of the key planning problems faced by railway operators. Crew costs typically constitute a substantial share of total expenditures: at Netherlands Railways, personnel costs account for approximately half of total operating expenses \citep{ns2025annual}. The problem is well-studied in the operations research literature, and we refer to \citet{heil2020railway} for a comprehensive review. Decision support tools based on optimisation models are used by major European operators, including Deutsche Bahn \citep{jutte2011optimizing} and Netherlands Railways \citep{abbink2005reinventing,huisman2024operations}. 

We study a canonical version of railway crew scheduling. Given a set of tasks, i.e., scheduled train trips, the goal is to assign them to crew members in the form of duties, sequences of tasks constituting a day of work, at minimum cost, while satisfying labour regulations. In particular, we require that every duty starts and ends at the same depot, does not exceed a maximum length, and contains a meal break within a bounded time of its start and end. These three rules are representative of those encountered in practice.

Column generation is the dominant solution method for railway crew scheduling \citep{heil2020railway}. The problem is naturally formulated as a large-scale set covering problem, where every column corresponds to a duty. Column generation solves the linear programming relaxation by iteratively generating duties of negative reduced cost. Embedding column generation within a branch-and-bound scheme yields an exact branch-and-price algorithm. It is common, however, to combine column generation with heuristic fixing strategies \citep{abbink2011solving, jutte2011optimizing, breugem2022column, dollevoet2024column}, yielding high-quality primal solutions, but no mechanism for improving the dual bound.

Exact methods, by contrast, remain scarce, and rarely solve instances containing more than a few hundred tasks \citep{heil2020railway}. \citet{boschetti2004exact} address the multi-depot case using a set partitioning formulation, iteratively enumerating duties of low reduced cost with respect to a dual bound and solving a reduced problem until optimality is certified. Other exact approaches rely on compact, flow-based formulations. \citet{fischetti2001polyhedral} develop a branch-and-cut approach for a single-commodity flow formulation of a simplified single-depot crew scheduling problem. \citet{hoffmann2019valid} propose valid inequalities for a crew-indexed multi-commodity flow formulation. However, these models struggle to capture complex duty rules, which are more naturally handled by the column-based approach. Such rules must instead be imposed through side constraints, leading to weak linear programming relaxations.

The status of exact methods contrasts sharply with vehicle routing, which shares considerable structure with railway crew scheduling: both are naturally formulated as set covering or partitioning problems and solved by column generation, with pricing modelled as a resource-constrained shortest path problem. The vehicle routing literature contains a large body of work on exact branch-price-and-cut methods, routinely solving instances with hundreds of customers to optimality \citep{costa2019exact}. Despite these structural similarities, the advances in vehicle routing have not yet translated to railway crew scheduling.

We believe two factors have contributed to this gap. First, the relative contribution of branch-and-price components, such as branching rules, primal heuristics, and reduced cost fixing, to solution quality and computational performance has not been systematically studied. Second, publicly available benchmark instances are scarce: real operator data is commercially sensitive, and labour regulations can vary across operators. To date, only a handful of data sets have been made available \Citep{neufeld2021efficient,rossum2022fast}, hindering a systematic comparison of solution methods.

This paper addresses both issues. Our first contribution is a systematic computational study of branch and price for railway crew scheduling, aimed at establishing which components of modern exact methods matter for this problem. To this end, we develop a branch-and-price algorithm that includes crew scheduling counterparts of the main components of state-of-the-art vehicle routing solvers, namely strong branching, primal heuristics, and reduced cost fixing, complemented by tailored branching rules and pricing. We evaluate the impact of each component, which to our knowledge has not been done before for railway crew scheduling. Moreover, we benchmark the algorithm against a multi-commodity flow formulation strengthened with symmetry-breaking constraints and a new family of valid inequalities.

Our second contribution is a public benchmark library of 100 instances to support our computational study, generated by mimicking the steps of a typical railway planning process, with every instance available in both a single-depot and a multi-depot variant. The library spans a wide range of instance sizes, is validated against historical planning data from Netherlands Railways, and is available with best known solutions and a solution checker at \url{https://github.com/BartvanRossum/rcsp-benchmark}. The instances can be used for benchmarking new solution methods or as a basis for studying new problem variants. 

Our experiments yield several findings, and point to a fundamental difference between exact methods for railway crew scheduling and vehicle routing. Branch and price returns high-quality solutions with average gaps below 1\% across all instance sizes, whereas the multi-commodity flow formulation is not competitive, failing to find a feasible solution on many of the smallest instances. However, these small gaps are driven primarily by the quality of primal solutions rather than by improvement in the dual bound, providing computational support for the widespread use of column generation-based primal heuristics. Semi-strong branching performs well but contributes little to closing the optimality gap, particularly on larger instances. The computational bottleneck is the restricted master problem rather than pricing, in contrast to vehicle routing, and reduced cost fixing has negligible impact. Together, these findings suggest that exact methods for railway crew scheduling face a different bottleneck than in vehicle routing: improving the dual bound appears to be the main challenge.

The remainder of this paper is structured as follows. In Section~\ref{sec:problem}, we formally describe the railway crew scheduling problem and present a set covering formulation. In Section~\ref{sec:method}, we present the branch-and-price algorithm and its components. In Section~\ref{sec:synthetic}, we introduce the benchmark library and validate it against real-life data from Netherlands Railways. In Section~\ref{sec:experiments}, we report the results of our computational experiments. We conclude and identify directions for future research in Section~\ref{sec:conclusion}.

\section{Problem Description}
\label{sec:problem}

Railway crew scheduling is defined on a directed acyclic graph $G = (V, A)$, where $V = \{s, t\} \cup K$ and $A$ is a set of arcs or \emph{connections}. Here, $K$ denotes the set of tasks to be covered, and $s$ and $t$ are source and sink nodes representing the depots at the start and end of the day, respectively. Let $D$ denote the set of stations that may serve as a crew base, with the single-depot case $|D|=1$ being a special instance of the formulation below. Each task $k \in K$ corresponds to a scheduled train trip and has a fixed start time, end time, start location, and end location, where the end time strictly exceeds the start time. The arc set $A$ contains connections between tasks that can be performed in order, where a necessary condition for an arc $(i, j) \in A$ is that task $j$ starts at the end location of task $i$ and that $j$ starts no earlier than $i$ ends. Additional restrictions may be imposed on the set of connections, such as a minimum transfer time when the rolling stock changes between tasks. For every crew base $d \in D$, arcs $(s, k)$ and $(k, t)$ indicate that task $k \in K$ is the first or last task of a duty starting or ending at $d$, respectively, and exist only for tasks starting or ending at that base. Since every task has a strictly positive duration, the graph is acyclic.

A \emph{duty} is an $s,t$-path in $G$, representing a sequence of tasks performed by a single crew member over the course of a working day. Not every $s,t$-path constitutes a valid duty: each duty must satisfy a set of labour rules. We consider three rules representative of those encountered in practice. First, every duty must start and end at the same crew base $d \in D$. Second, the total length of a duty, defined as the time between the start of the first task and the end of the last task, may not exceed a prescribed maximum duty length $L$. Third, every duty must contain a meal break, and the time elapsed before and after the break may not exceed a prescribed maximum time without break $L^B$. 

To model the meal break constraint, we designate a subset of arcs $A^B \subseteq A$ as \emph{break arcs}: a break arc indicates that a meal break can be taken during the corresponding connection, and specifies the time at which the break starts and ends. A duty is feasible with respect to the meal break constraint if it contains a break arc for which the time from the start of the duty to the beginning of the break, and from the end of the break to the end of the duty, both do not exceed the maximum time without break.

The cost of a duty is the sum of the costs of its connections. In practice, this typically consists of a variable component proportional to duty length and a fixed cost. The railway crew scheduling problem is to select a minimum-cost set of duties that cover all tasks at least once.

\subsection{Set Covering Formulation}
\label{subsec:problem_formulation}

Let $\Delta$ denote the set of all feasible duties, and let $c_\delta \geq 0$ denote the cost of duty $\delta \in \Delta$. For each task $k \in K$ and duty $\delta \in \Delta$, let $a_{\delta k} \in \{0,1\}$ indicate whether task $k$ is covered by duty $\delta$. Introducing a binary decision variable $x_\delta$ to indicate whether duty $\delta$ is selected, the railway crew scheduling problem can be formulated as
\begin{subequations}
\begin{align}
\min \quad & \sum_{\delta \in \Delta} c_\delta x_\delta & \label{eq:objective} \\ 
\text{s.t.} \quad & \sum_{\delta \in \Delta} a_{\delta k} x_\delta \geq 1 && \forall k \in K \label{eq:cover} \\ 
& x_{\delta} \in \{0, 1\} && \forall \delta \in \Delta. \label{eq:domain}
\end{align}
\label{eq:set_cover}%
\end{subequations}
The objective~\eqref{eq:objective} is to minimise total duty costs. Constraints~\eqref{eq:cover} enforce that each task is covered by at least one duty, and the binary domain is specified by~\eqref{eq:domain}. The formulation allows overcovering: a task may appear in more than one selected duty, corresponding to a crew member \emph{deadheading}, i.e., moving through the network as a passenger on a trip they are not operating themselves. This differs from airline crew pairing, where overcovering is typically penalised or forbidden. Since the number of feasible duties $|\Delta|$ grows exponentially with instance size, the linear programming (LP) relaxation of~\eqref{eq:set_cover} is solved via column generation, and optimal integer solutions can be obtained using branch and price.

\section{Branch and Price}
\label{sec:method}

Branch and price is a branch-and-bound algorithm in which the LP relaxation at every node of the search tree is solved via column generation. It is the exact method of choice for set covering and set partitioning problems where the number of variables is too large to handle explicitly, as is the case in railway crew scheduling. We discuss the components of our branch-and-price algorithm in turn: the column generation routine in Section~\ref{subsec:method_colgen}, branching rules in Section~\ref{subsec:method_branching}, primal heuristics in Section~\ref{subsec:method_primal}, and reduced cost-based variable fixing in Section~\ref{subsec:method_fixing}. We refer the interested reader to \citet{dollevoet2024column} for a detailed introduction to column generation for railway crew scheduling, and to \citet{desrosiers2024branch} for more details on branch and price.

\subsection{Column Generation}
\label{subsec:method_colgen}

We solve the LP relaxation of~\eqref{eq:set_cover} using column generation, initialising the model with only a small subset of duties $\Delta' \subseteq \Delta$, yielding the restricted master problem (RMP):
\begin{subequations}
\begin{align}
\min \quad & \sum_{\delta \in \Delta'} c_\delta x_\delta & \label{eq:rmp_objective} \\ 
\text{s.t.} \quad & \sum_{\delta \in \Delta'} a_{\delta k} x_\delta \geq 1 && \forall k \in K \label{eq:rmp_cover} \\ 
& x_{\delta} \geq 0 && \forall \delta \in \Delta'. \label{eq:rmp_domain}
\end{align}
\label{eq:RMP}%
\end{subequations}
Column generation iterates between solving the RMP and a pricing problem, whose goal is to identify a duty with negative reduced cost. Let $\bm{\lambda} \in \mathbb{R}^{|K|}$ denote the vector of dual variables associated with constraints~\eqref{eq:rmp_cover}. The reduced cost of duty $\delta \in \Delta$ is given by
\begin{align}
    \tilde{c}_{\delta} = c_{\delta} - \sum_{k \in K} a_{\delta k} \lambda_{k}.
    \label{eq:reduced_cost}
\end{align}
If a negative reduced cost duty exists, it is added to $\Delta'$ and the RMP is re-solved. If not, the current RMP solution is provably optimal for the LP relaxation of~\eqref{eq:set_cover}.

\subsubsection{Pricing Algorithm}
\label{subsubsec:method_pricing}

Observe that the reduced cost~\eqref{eq:reduced_cost} decomposes over the connections of a duty: $c_\delta$ is a sum of connection costs, and the dual term $-\lambda_k$ can be placed on all incoming connections of task $k$. Due to the specific structure of the duty rules and the graph $G$, the pricing problem admits tailored polynomial-time algorithms \Citep{huisman2007column,rossum2022fast}, exploiting the fact that shortest paths on $G$ can be determined in time linear in the number of arcs. 

However, we opt to model the pricing problem as the problem of finding a resource-constrained shortest $s,t$-path in $G$, where the resources enforce the labour rules described in Section~\ref{sec:problem}, and solve it with a labelling algorithm. This has the following benefits. First, it allows us to easily generate multiple negative reduced cost columns per iteration, which accelerates column generation. Second, labelling algorithms admit a wide range of acceleration techniques, such as the use of completion bounds. Third, the labelling framework can be reused for reduced cost-based duty enumeration and arc fixing, as described in Section~\ref{subsec:method_fixing}. Since the pricing algorithm does not dominate runtime on the instances under consideration, we do not use heuristic pricing in our experiments.

In the multi-depot case, every duty must still start and end at the same crew base, even though any station in $D$ may serve as one. Instead of tracking the crew base as an additional resource in the label, we enforce this by decomposing the pricing problem by crew base: for a given base $d \in D$, we solve a separate pricing problem in which the source and sink arcs are restricted to $d$, so that any $s,t$-path necessarily starts and ends at $d$. 

We use the bidirectional labelling algorithm of \citet{righini2006symmetry}, where labels are extended in two directions starting at both the source and the sink. Every label corresponds to a partial duty and tracks four quantities: the accumulated reduced cost, the elapsed duty length, the time elapsed since the last meal break, and an indicator of whether a meal break has been taken. The duty length and the time since the last break jointly determine the position of the most recent break, if any. A label $\ell$ dominates a label $\ell'$ ending at the same task if $\ell$ has lower or equal reduced cost, lower or equal duty length, lower or equal time without break, and has taken a meal break whenever $\ell'$ has. In this case, label $\ell'$ can be discarded without loss of optimality. We refer to Appendix~\ref{app:labelling} for a detailed description of the label definition, extension, and dominance check.

We apply various techniques to accelerate the labelling algorithm. First, we compute completion bounds via simple shortest path problems in $G$, providing lower bounds on the reduced cost, time without break, and duty length achievable from any task. The first completion bound must be computed in every column generation iteration, whereas the latter two can be pre-computed before the start of the column generation algorithm. Labels whose best-case completion attains nonnegative reduced cost or violates a resource constraint are pruned. Second, we use dynamic midpoint selection \citep{tilk2017asymmetry} to balance the number of generated forward and backward labels. Third, we decompose the pricing problem by starting hour: for each hour of the day, we solve a separate pricing problem in which source arcs are restricted to tasks starting within that hour. In the multi-depot case, this decomposition is applied within each crew base's pricing problem, so that subproblems are indexed by both crew base and starting hour. These subproblems are solved in parallel.

\subsubsection{Acceleration Techniques}
\label{subsubsec:method_accelerate}

We apply three well-known acceleration techniques to the column generation procedure. First, we use Neame's dual smoothing \citep{pessoa2018automation}, which stabilises the dual variables across iterations and accelerates overall convergence. Second, we apply a column selection strategy in which generated columns are sorted by increasing reduced cost, and only columns that are sufficiently task-disjoint from previously selected ones are added to the RMP \citep{breugem2022column}. This ensures that the columns added in each iteration cover a diverse set of tasks and avoids adding redundant columns to the RMP. Third, we perform column management by periodically removing columns from the RMP that have not been in the basis for a given number of iterations, keeping at least a certain number of non-basic columns. 

\subsection{Branching}
\label{subsec:method_branching}

At every fractional node of the branch-and-price tree, we apply one of the following branching rules to obtain a set of branching candidates:
\begin{enumerate}
    \item \emph{Number of duties.} Branch on the total number of duties, i.e., $\sum_{\delta \in \Delta'} x_\delta$.
    \item \emph{Connections.} Branch on the fractional value of a single connection, i.e., the total flow through arc $(i,j) \in A$ across all selected duties.
    \item \emph{Connection time.} For a given task, sort all incoming or outgoing connections by start or end time and branch on whether the total flow through connections before or after a given time threshold exceeds an integer value.
    \item \emph{Depot connections.} Same as rule (2), but restricted to connections to or from the depot.
    \item \emph{Cumulative duties started.} Branch on the cumulative number of duties that have started at or before a given time of day.
    \item \emph{Station flow.} Branch on the total flow of connections arriving at or departing from a given station within a specified time window.
    \item \emph{Depot flow.} Same as rule (6), but restricted to the depot.
\end{enumerate}
Rules (3), (5), (6), and (7) branch on aggregates of connections rather than individual ones, with the aim of producing a more balanced search tree. Rule (4) focuses on depot connections, since solution cost is uniquely determined by the source and sink arcs.

A single branching rule typically yields multiple candidates. We consider the following candidate selection rules:
\begin{enumerate}
    \item \emph{Most fractional.} Select the candidate whose fractional value is closest to $0.5$.
    \item \emph{Semi-strong branching.} Evaluate both branches of a candidate by solving the RMP with the branching constraint imposed, without performing additional pricing iterations, and branch on the candidate yielding the highest weighted lower bound improvement.
    \item \emph{Full strong branching.} A restricted number of candidates is first identified using semi-strong branching, after which a limited number of column generation iterations are performed on each candidate to sharpen the bound estimate.
\end{enumerate}
The goal of strong branching is to avoid ineffective branching decisions by investing some computational effort to identify the most promising branching candidate. While this may increase the time per node, it may still reduce overall computation time by yielding a smaller search tree \citep{achterberg2005branching}. For both strong branching variants, the score of a candidate is computed as $\frac{1}{6} \Delta_{\max} + \frac{5}{6} \Delta_{\min}$, where $\Delta_{\max}$ and $\Delta_{\min}$ denote the larger and smaller bound improvements on the two sides, respectively, and the candidate with the highest score is selected. In semi-strong branching, multiple candidates are evaluated in parallel. 

\subsection{Primal Heuristics}
\label{subsec:method_primal}

We apply primal heuristics at the root node to obtain high-quality upper bounds early on. These heuristics frequently return the best solutions found during the entire branch-and-price algorithm. We consider two classes of primal heuristics: restricted master heuristics and diving heuristics.

The \emph{restricted master heuristic} (RMH) first solves the LP relaxation via column generation, and then solves the IP restricted to the columns currently in the RMP. The solution quality depends on whether the RMP contains a sufficiently rich set of columns. To improve this, we can augment the RMP with additional columns of low reduced cost prior to solving the IP (see Section~\ref{subsec:method_fixing}).

An effective alternative is \emph{diving}, a column generation-based fixing heuristic that is frequently used in crew scheduling \citep{heil2020railway}. Diving iteratively selects one or more fractional duties, fixes them to one, and re-solves the LP via column generation, repeating until an integer solution is found. Due to the set covering structure of~\eqref{eq:set_cover}, this procedure is guaranteed to find a feasible solution. This does not hold for set partitioning formulations as used in vehicle routing.

A more sophisticated diving algorithm uses \emph{backtracking} \citep{sadykov2019primal}, exploring a larger portion of the diving tree at the cost of additional computation time. The scheme is governed by two parameters: $\texttt{maxDepth}$ and $\texttt{maxDiscrepancy}$. Backtracking is performed up to a depth of $\texttt{maxDepth}$ in the diving tree: when a branch is backtracked, the duty that was fixed in that branch is added to a tabu list and forbidden from being selected in any sibling or descendant branch. Backtracking continues as long as the size of the tabu list does not exceed $\texttt{maxDiscrepancy}$. Together, these two parameters control the width of the search: $\texttt{maxDepth}$ limits how far back the algorithm can retrace its steps, while $\texttt{maxDiscrepancy}$ bounds the total number of alternative branching decisions explored. 

Finally, we use a diving-based \emph{postprocessing heuristic} to improve the best solution returned by branch and price. Given a feasible solution, a fraction of its duties is fixed to one and the remaining tasks are re-optimised using diving. This approach is essentially a form of local search. We consider three rules for selecting which duties to fix. 
\begin{enumerate}
    \item \emph{Random.} Fix a randomly chosen subset of duties.
    \item \emph{Cost-per-task ratio.} Select duties with the \emph{lowest cost-per-task ratio}, keeping the most efficient duties fixed and allowing the relatively expensive ones to be replaced.
    \item \emph{Time window.} Randomly select a reference duty and fix the duties that are furthest from it in time, with the aim of increasing exchange opportunities among unfixed duties.
\end{enumerate}
The postprocessing heuristic is applied iteratively for a fixed number of iterations, with each iteration starting from the best solution found so far.

\subsection{Reduced Cost Fixing}
\label{subsec:method_fixing}

Reduced cost fixing is a well-known technique to accelerate branch-and-bound algorithms \citep{crowder1983solving}, and is often employed in vehicle routing \citep{costa2019exact}. The idea is to use the gap between the best known primal solution and the LP relaxation at the current node to fix variables to zero. Formally, let $z^*$ denote the objective value of the best known primal solution and let $z_{\text{LP}}$ denote the LP relaxation value at the current node. Any duty $\delta \in \Delta$ with reduced cost $\tilde{c}_\delta \geq z^* - z_{\text{LP}}$ can be fixed to zero without loss of optimality, since including it in any solution would result in an objective value of at least $z^*$. We apply reduced cost fixing at two levels: entire duties and individual connections.

\emph{Duty enumeration.} Once the gap is sufficiently small, we enumerate all duties with reduced cost below $z^* - z_{\text{LP}}$ using a modified version of the bidirectional labelling algorithm. We first run the exact labelling algorithm with dominance to identify which connections are used by at least one duty with sufficiently low reduced cost. We then perform an exhaustive labelling without dominance, restricted to the active connections identified in the first pass, to enumerate all such duties. The enumerated duties are added to the RMP, which is then solved as a regular IP, replacing further branching at this node. If the lower bound of this program is at least $z^*$, the node can be pruned.

\emph{Arc fixing.} When the gap is not small enough for duty enumeration to be tractable, we can attempt to fix individual connections. For each connection $(i,j) \in A$, we compute the minimum reduced cost of any duty passing through it using a modified version of the bidirectional labelling algorithm. If this minimum reduced cost exceeds $z^* - z_{\text{LP}}$, the connection can be removed from $G$ without loss of optimality, thereby accelerating subsequent pricing iterations. 

\section{Synthetic Benchmark Instances}
\label{sec:synthetic}

We introduce a library of benchmark instances for railway crew scheduling, designed to reflect the structure of real-world instances. Every instance is provided in a single-depot and a multi-depot variant, which share the same tasks and connections between tasks but differ in the set of stations that may serve as a crew base. We describe the generation procedure in Section~\ref{subsec:synthetic_generation}, report aggregate statistics of the resulting instances in Section~\ref{subsec:synthetic_statistics}, and validate them against real-life data from Netherlands Railways in Section~\ref{subsec:synthetic_validation}.

\subsection{Generation Procedure}
\label{subsec:synthetic_generation}

We mimic the planning process that a railway operator typically follows, including line planning, timetabling, and rolling stock scheduling \citep{huisman2005operations,borndorfer2018handbook}. Each step takes the output of the previous step as input, ensuring that the resulting instances are constructed in the same way as real-life crew scheduling instances. Figure~\ref{fig:generation} illustrates the generation pipeline.

\input{figs/generation.tex}

The procedure begins by generating a physical infrastructure network of stations connected by tracks. Based on station size and pairwise distance, passenger demand data are then derived using a gravity model, reflecting the number of passengers wishing to travel between each pair of stations. Given this demand, a minimum-cost line plan is constructed, specifying the frequencies and itineraries of trains operating in the network. A periodic timetabling problem is then solved to assign arrival and departure times to all train trips, minimising overall transfer times. The resulting timetable is rolled out over the full operating day and a rolling stock schedule is generated greedily, which determines where rolling stock changes occur between consecutive tasks. From the rolling stock schedule, tasks and connections are extracted, forming the necessary input for crew scheduling. A postprocessing step removes tasks that cannot be covered in a feasible solution. We describe each step in detail below.

\paragraph{Infrastructure network.} 
We randomly sample stations on a $500 \times 500$ grid, where one distance unit corresponds to one minute of travel time. Large stations serve as the start and end points of lines and as crew transfer locations, whereas small stations are intermediate stops only. Large stations are placed by repeatedly sampling uniformly random coordinates until each new station is at least 10 and at most 50 distance units from at least one previously placed station. The track network is constructed in two steps: (i) compute the minimum spanning tree over the large stations, where track travel times equal the Euclidean distance in minutes, and (ii) greedily add shortcut arcs, where each arc is chosen to maximally reduce total pairwise travel time through the network and no two arcs may cross. Small stations are then placed along existing tracks: a random coordinate is sampled and projected onto the nearest track segment, subject to a separation of at least 5 and at most 15 units from the nearest station.

\paragraph{Passenger demand.}
Each large station is assigned a population drawn uniformly from $[0.5, 1.0]$, and each small station from $[0.05, 0.25]$. Origin-destination (OD) demand is generated using a simple gravity model: demand between two large stations is proportional to $p_o \cdot p_d / d_{od}$, where $p_o$ and $p_d$ are the populations of the origin and destination, and $d_{od}$ is their Euclidean distance. Small stations only generate demand towards their two nearest large stations. OD pairs with a travel distance exceeding 120 minutes are excluded, demand is normalised to equal the sum of station populations, and OD pairs with normalised demand below 0.0025 are discarded.

\paragraph{Line planning.}
A line pool is constructed by enumerating all paths in the infrastructure network between pairs of large stations with a travel time of at most 210 minutes, a detour factor of at most 2.5 relative to the shortest path, and containing four or more stations. At most two lines are retained per station pair, preferring lines that pass through more large stations. Given the line pool, a change-and-go network (CGN) is constructed to model how passengers route through the network using the available lines \citep{schobel2012line}. We impose a transfer penalty of 10 minutes for changing lines. For each OD pair, passenger routes are enumerated using depth-first search in the CGN, starting from the shortest path and incrementally relaxing the detour factor in steps of 0.05 up to a maximum of 1.5, retaining at most 100 routes per OD pair sorted by travel time. We select lines by solving a line planning model, where the goal is to minimise line costs subject to coverage of all passenger demand and capacity constraints on each line segment. The output is a set of selected lines with associated frequencies, which forms the input to the timetabling step.

\paragraph{Timetabling.}
A periodic, hourly timetable is computed using the Periodic Event Scheduling Problem (PESP) formulation \citep{liebchen2007timetable}, minimising a weighted sum of transfer times. The resulting periodic timetable is expanded into a full operating day by replicating each line departure at its scheduled frequency between 06:00 and 22:00, yielding a set of scheduled trips with fixed departure and arrival times at each station.

\paragraph{Rolling stock scheduling.}
Rolling stock units are assigned to scheduled trips using a greedy procedure. Each unit starts at the departure station of a line and is matched to return trips on the reverse line subject to a minimum turnaround time. If a unit cannot return to its starting station after its last trip, a deadhead return trip is added. The output is a set of trips with known rolling stock assignments, which determines where rolling stock changes occur and thereby affects connection feasibility in the crew task network.

\paragraph{Tasks and connections.}
Tasks correspond to segments between consecutive large stations on a scheduled trip. For the single-depot variant, one of the large stations is designated as the crew base. We select the station with the highest degree in the infrastructure network, reflecting the practice of locating crew bases at major junctions. Source and sink connections are only created for tasks starting and ending at the crew base, respectively. A connection between two tasks at the same station is feasible if the transfer time is at most 2 hours. If the rolling stock changes between tasks, an additional minimum transfer time of 10 minutes is imposed. A connection is marked as a break arc if the transfer gap is at least 40 minutes, allowing a 30-minute meal break with 5 minutes of travel to the station canteen on each side. Connection costs consist of a variable component of 1 unit per second of duty time, plus a fixed cost equal to 8 hours of variable cost charged on source connections.

\paragraph{Postprocessing.}
After constructing the task network, the LP relaxation of the crew scheduling problem is solved via column generation. Any task that cannot be covered is removed and the task network is rebuilt. This step ensures that every instance in the benchmark admits a feasible solution. Tasks that are isolated in the task network, i.e., that have no incoming or outgoing connections, are also removed iteratively until no such tasks remain. Typically, no more than a handful of tasks are removed in this step.

\paragraph{Multi-depot instances.} For every single-depot instance, we construct a multi-depot counterpart in which every large station may serve as a crew base, rather than only the one selected for the single-depot variant. Tasks and connections between tasks are unchanged: we only add, for every large station that is not already the single-depot crew base, the source and sink connections that would have been generated for it had it been selected as the depot, using the same cost and duration rules described above. Every connection in the single-depot instance therefore also appears in its multi-depot counterpart, and for a given instance, its single-depot solution remains feasible in the multi-depot case.

\subsection{Library}
\label{subsec:synthetic_statistics}

We obtain a benchmark library of 100 unique instances by applying the instance generation procedure with a varying number of large stations and different random seeds, each provided in a single-depot and a multi-depot variant. The instances are partitioned into four groups based on the number of tasks: tiny (0-500 tasks), small (501-1,000 tasks), medium (1,001-2,000 tasks), and large (more than 2,000 tasks), with 25 instances per group. This partitioning is used consistently throughout all computational experiments. Table~\ref{tab:instances} reports aggregate statistics per group, with the average number of connections given separately for the single-depot and multi-depot variants, and Figure~\ref{fig:instances} shows a selection of synthetic infrastructure networks. The instances are publicly available, together with best known solutions and a solution checker, at \url{https://github.com/BartvanRossum/rcsp-benchmark}.

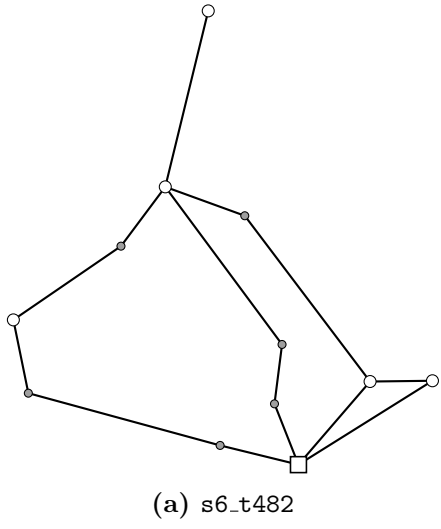
\input{tables/instanceTable}

\input{figs/instancePlot}

\subsection{Validation against Real-Life Data}
\label{subsec:synthetic_validation}

To assess whether the synthetic instances are representative of real-life crew scheduling instances, we compare them against instances derived from historical crew schedules of Netherlands Railways (NS). For each day of a representative week in October 2021, we sample random combinations of 1 to 6 crew bases, yielding 42 NS instances in total, and construct crew scheduling instances using the duties assigned to these crew bases. To convert them to single-depot instances, we designate the crew base with the highest number of duties as the depot. For each instance, we solve the root node LP relaxation via column generation and record six metrics: average task duration, average connection degree, root node lower bound per task, average number of duties in the LP solution, average duty length, and root node solve time.

\input{figs/validation}

Figure~\ref{fig:validation} compares these metrics across synthetic and NS instances as a function of instance size, as measured by the number of tasks. The distributions match well overall: average task duration, lower bound per task, average number of duties, and average duty length are closely aligned across the two sets. Two metrics show a modest discrepancy: the synthetic instances have a slightly higher connectivity than the NS instances of similar size, and a somewhat higher root node solve time. The former is likely a consequence of the synthetic network topology, whereas the latter is likely a result of the higher connectivity. Despite these differences, the general patterns are consistent, and the synthetic instances are computationally at least as challenging as the NS instances. We therefore consider them a suitable basis for our computational experiments.

\section{Computational Experiments}
\label{sec:experiments}

We evaluate the performance of the branch-and-price algorithm and its individual components on the synthetic benchmark instances introduced in Section~\ref{sec:synthetic}. Sections~\ref{sec:experiments_setup}-\ref{subsec:experiments_mcf} focus on the single-depot instances, and we apply the resulting best-performing configuration to the multi-depot instances in Section~\ref{subsec:experiments_multi_depot}. We describe the experimental set-up in Section~\ref{sec:experiments_setup}, and analyse the effect of the branching rule and candidate selection strategy in Section~\ref{subsec:experiments_branching}. We evaluate the primal heuristics in Section~\ref{subsec:experiments_primal}, and report full branch-and-price runs including an assessment of reduced cost fixing in Section~\ref{subsec:experiments_branchPrice}. We then evaluate the postprocessing heuristic in Section~\ref{subsec:experiments_postprocess}, and compare branch and price against a multi-commodity flow formulation in Section~\ref{subsec:experiments_mcf}. Finally, we evaluate the algorithm on the multi-depot instances in Section~\ref{subsec:experiments_multi_depot}, and report best known bounds and solutions for all instances in Appendix~\ref{app:bestKnown}.

\subsection{Set-Up}
\label{sec:experiments_setup}

All experiments use the following crew scheduling parameters. The maximum duty length is set to 9.5 hours, the meal break length to 30 minutes, and the maximum time without a meal break to 5.5 hours. Connection costs consist of a variable component of one unit per second, plus a fixed activation cost equal to 8 hours of variable cost. These parameters are representative of those used in practice by Netherlands Railways.

The column generation procedure uses dual smoothing with self-tuning smoothing parameter $\alpha$, and a task-disjoint column selection strategy in which every newly added column must cover at least 50\% of its tasks uniquely relative to all previously selected columns in the same iteration. Every 5 iterations, all columns that have not been in the basis for the last 10 iterations are removed, keeping at least 1,000 non-basic columns at all times. These parameters are based on the parameter sensitivity analysis in Appendix~\ref{sec:parameter}. For semi-strong branching, up to 10 candidates are evaluated. For full strong branching, the 10 semi-strong candidates are reduced to 3, each of which is then evaluated using 5 column generation iterations per branch.

We implemented all algorithms in \texttt{Java}. All LPs are solved using \texttt{CPLEX 22.1.0} on six threads. Pricing subproblems are solved in parallel on eight threads. All experiments are conducted on cluster nodes with 16 GB RAM and an AMD Rome 7H12 processor.

\subsection{Branching Rules}
\label{subsec:experiments_branching}

We test all combinations of branching rule and candidate selection strategy, using a time limit of one hour. Table~\ref{tab:branching} reports the average lower bound improvement over the root node ($\Delta$LB) and the number of nodes processed for each combination, averaged over instances of similar size. Every setting includes branching on duties and connections as the base rules, the additional branching rule is indicated in the table. For branching rules that involve a time window, i.e., branching rules (5)-(7), we use a window of 15 minutes.

\input{tables/branchingTable}

The results show a clear advantage for semi-strong branching over both most fractional and full strong branching. Without strong branching, the average $\Delta$LB on small instances reaches at most 0.151\%, while semi-strong branching achieves 0.258\% with the default branching scheme. Full strong branching achieves a comparable $\Delta$LB of 0.232\% on small instances, but processes far fewer nodes: 263 compared to 642 for semi-strong branching. This indicates that full strong branching requires too much time per node, leaving insufficient time to explore the search tree. The same pattern holds on the medium and large instances.

None of the additional rules consistently improves over the default scheme of branching on duties and connections. Under semi-strong branching, the default scheme achieves the highest $\Delta$LB of 0.258\%, 0.089\%, and 0.033\% on small, medium, and large instances respectively. The connection time rule is the most competitive alternative, matching the default on medium instances but falling short on the small and large ones. On the tiny instances, connection time and station flow marginally exceed the default, at 0.600\% against 0.599\%, but the difference is too small to be meaningful. The remaining rules, depot connection, cumulative start, and depot flow, reduce $\Delta$LB across all instance sizes, suggesting that they introduce branching decisions that are less effective than those produced by the default scheme.

Two broader observations are worth noting. First, the effect of branching diminishes as instances grow: under semi-strong branching, $\Delta$LB decreases from 0.599\% on tiny instances to 0.258\% on small, 0.089\% on medium, and 0.033\% on large instances, indicating that branching contributes little to closing the optimality gap on larger instances. Second, the relative rankings across configurations are stable across instance sizes. Based on these results, we use semi-strong branching with the default branching scheme in all subsequent experiments.

\input{figs/connectionRemovalScatter}

To investigate why branching performs so poorly at improving the dual bound, we examine the degree to which the root LP solution is degenerate, i.e., how many alternative optimal solutions exist. To this end, we consider subsets of up to two connections used in the root LP solution, forbid all duties using these connections, and re-solve the restricted master problem without further column generation. We count how many optimal solutions that differ in the set of connections used can be obtained in this way, across all subsets explored within a one-hour time limit. As the number of RMP re-solves grows rapidly, we only consider the 25 tiny instances. Since we do not re-solve using column generation, impose a time limit, and only consider subsets of size at most two, the resulting count is a lower bound on the true number of such optimal solutions.

Figure~\ref{fig:connection_removal_scatter} reports the results. For the sake of clarity, we have omitted from the plot one instance with 333 tasks where the count already exceeded 22 thousand. The number of alternative optimal solutions found appears to grow exponentially with the number of tasks. On average, 84\% of the connections used are common to every solution found for a given instance, while the remaining 16\% vary across solutions. As the LP relaxation admits a large number of alternative optima, branching on individual connections is unlikely to meaningfully improve the dual bound. In particular, forcing one connection to an integer value still leaves many equally good alternatives available.

\subsection{Primal Heuristics}
\label{subsec:experiments_primal}

We evaluate all primal heuristics with the following set-up. The RMH is given a time limit of 30 minutes, and is optionally augmented with duty enumeration using columns with reduced cost up to 500. Diving variants use combinations of (\texttt{maxDepth}, \texttt{maxDiscrepancy}) equal to (0, 0), (1, 1), (1, 2), and (2, 1), where the first option is default diving and the last three options include backtracking. Table~\ref{tab:primal} reports the optimality gap relative to the root node lower bound and the computation time for each variant, averaged over instances of similar size.

\input{tables/primalTable}

Diving outperforms the RMH across all instance sizes except the tiny ones, where RMH with enumeration performs best. The RMH achieves gaps of 1.24\%, 2.86\%, 4.55\%, and 6.57\% on tiny, small, medium, and large instances respectively, hitting the time limit on all but the tiny and small instances. Enumerating and adding duties improves the gap on tiny and small instances, from 1.24\% to 0.88\% and from 2.86\% to 1.77\% respectively, but worsens it on medium and large instances, where the larger column pool increases the difficulty of solving the restricted IP and actually makes it harder to find good primal solutions. Default diving finds solutions with gaps between 1.06\% and 1.46\% in 4 to 584 seconds depending on instance size, requiring far less computation time than the RMH.

Backtracking improves over default diving at a modest additional computational cost. Backtracking (1, 1) reduces the gap by 0.34, 0.27, 0.03, and 0.08 percentage points on tiny, small, medium, and large instances respectively, to 1.12\%, 1.07\%, 1.03\%, and 0.98\%. Increasing the discrepancy parameter from 1 to 2, i.e., comparing (1, 1) to (1, 2), yields a further reduction on all but the tiny instances, to 0.92\%, 1.00\%, and 0.95\% on small, medium, and large instances, at roughly 50\% additional computation time. Increasing the depth parameter from 1 to 2, i.e., comparing (1, 1) to (2, 1), produces a similar improvement at comparable computation time, and performs slightly worse than (1, 2) on all but the large instances. Based on these results, we use backtracking (1, 2) as the primal heuristic in all subsequent branch-and-price experiments.

\subsection{Branch and Price}
\label{subsec:experiments_branchPrice}

We evaluate full branch-and-price runs using the components identified in the preceding sections: semi-strong branching with the default branching scheme, and backtracking (1, 2) as primal heuristic at the root node. In addition, every 20 nodes we call the restricted master heuristic with a time limit of one minute and without augmented columns. We compare two configurations: with and without reduced cost fixing. In the former, we apply arc fixing at every node where the gap is below 2,500 and attempt to close the node using duty enumeration whenever the gap is below 1,500 and the number of enumerated duties does not exceed 100,000. Both configurations use a time limit of three hours. Table~\ref{tab:solver_comparison} reports the average lower bound and upper bound, both normalised by the root node lower bound, the optimality gap, the number of nodes processed, and a breakdown of computation time across four components: RMP solving, pricing, branching, and calls to the primal heuristics at the root node and subsequent nodes. The remaining time is overhead, and is not reported in the table. The best known lower bounds across all experiments of Sections~\ref{subsec:experiments_branching} and~\ref{subsec:experiments_branchPrice} are reported in Tables~\ref{tab:best_tiny}-\ref{tab:best_large} in Appendix~\ref{app:bestKnown}.

\input{tables/tableSolver}

Overall gaps are below 1\% across all instance sizes, averaging 0.03\%, 0.58\%, 0.89\%, and 0.85\% on tiny, small, medium, and large instances respectively under the configuration without reduced cost fixing. These small gaps are driven primarily by the quality of the primal solutions rather than by improvement in the lower bound: the LB increase over the root node ranges from 0.59\% on tiny instances to just 0.03\% on large instances, confirming the observation from Section~\ref{subsec:experiments_branching} that branching contributes little to closing the gap.

The time breakdown shows that pricing accounts for less than 6\% of total computation time across all instance sizes, while RMP solving consumes 21-48\% and grows as a share of total time with instance size. Branching accounts for a further 23-39\%, largely explained by the use of semi-strong branching, and the primal heuristics for 23-59\%, mostly consisting of diving with backtracking at the root node. The dominance of the RMP is a direct consequence of the large number of columns accumulated over the course of the search, and stands in contrast to vehicle routing, where pricing is typically the bottleneck \citep{costa2019exact}.

Reduced cost fixing does not improve solution quality at any instance size: gaps differ by at most 0.08 percentage points between the two configurations, in either direction. Its only visible effect is on the tiny instances, where the gap is small enough to trigger duty enumeration regularly: the heuristic there accounts for 59.1\% of computation time versus 31.3\% without fixing, and the number of nodes processed drops from 606 to 367, without any improvement in the gap. The gap rarely falls below the thresholds required to trigger arc fixing or duty enumeration on the larger instances, explaining the absence of any effect at those sizes. This is confirmed by the best known lower and upper bounds reported in Appendix~\ref{app:bestKnown}: the remaining gaps are too large for reduced cost fixing to be applicable.

Table~\ref{tab:solver_comparison} also reports results on the 42 NS instances described in Section~\ref{subsec:synthetic_validation}. The average gap of 0.19\% without reduced cost fixing is in line with the synthetic instances of comparable size, consistent with the finding in Section~\ref{subsec:synthetic_validation} that the synthetic instances are representative of real-life crew scheduling data. As for the synthetic instances, reduced cost fixing has no discernible effect on the final gap.

\subsection{Postprocessing Heuristic}
\label{subsec:experiments_postprocess}

We evaluate the postprocessing heuristic described in Section~\ref{subsec:method_primal} on the best solutions found across the experiments of Sections~\ref{subsec:experiments_primal} and~\ref{subsec:experiments_branchPrice}. We consider all combinations of the three duty selection rules (random, cost-per-task, and time window) and three fixing fractions (25\%, 50\%, and 75\%), running each configuration for 25 iterations. Table~\ref{tab:local_search} reports the average improvement in upper bound ($\Delta$UB) and total computation time per configuration, averaged over instances of similar size. Since almost all tiny instances were already solved to optimality by the branch-and-price algorithm, we only report results for the small, medium, and large instances. Best known upper bounds after postprocessing are reported in Tables~\ref{tab:best_tiny}-\ref{tab:best_large} in Appendix~\ref{app:bestKnown}, and the corresponding solutions are available in the online repository.

\input{tables/tableLocalSearch}

The results show a clear and consistent pattern across all instance sizes: fixing a smaller fraction of duties yields larger improvements, at the cost of increased computation time. Fixing 25\% of duties with random selection achieves the largest improvements of 0.17\%, 0.18\%, and 0.16\% on small, medium, and large instances respectively, while fixing 75\% produces improvements of at most 0.04\% across all settings. This is not surprising: fixing more duties restricts the solution space, leaving less room for improvement.

Among the three selection rules, random selection generally outperforms the alternatives, followed closely by the time window rule. The cost-per-task rule performs substantially worse across all fixing fractions and instance sizes, achieving improvements up to three times smaller than random selection at 25\% fixing. This suggests that allowing for some degree of randomness is beneficial, compared to a fully deterministic fixing rule. The time window rule is competitive with random selection at 50\% and 75\% fixing, but falls short at 25\% fixing, where random selection has more freedom to identify improving moves.

\subsection{Comparison with Multi-Commodity Flow Formulation}
\label{subsec:experiments_mcf}

We conclude our computational study by comparing branch and price against an alternative exact method for railway crew scheduling. To this end, we consider the multi-commodity flow formulation presented in Appendix~\ref{app:mcf}, and solve it directly using \texttt{CPLEX 22.1.0}. We test two variants: one with the symmetry-breaking constraints only, and one that additionally includes the valid inequalities. The formulation is crew-indexed and therefore requires an a priori upper bound on the number of duties in an optimal solution. We set this bound to the number of duties in the best known solution of the instance at hand plus five. Such a bound would not be available to a practitioner, so this choice works in favour of the multi-commodity flow formulation. We compare against the branch-and-price algorithm of Section~\ref{subsec:experiments_branchPrice} without duty enumeration and reduced cost fixing, since these components were shown to have no discernible effect on solution quality. Both methods are given a time limit of three hours.

\input{tables/tableMCF}

Table~\ref{tab:mcf} reports the resulting bounds and optimality gaps for both methods, where the best gap and computation time per instance are indicated in bold. A dash indicates that no feasible solution was found within the time limit. Branch and price dominates both flow formulations. It solves 20 of the 25 tiny instances to optimality, compared to 1 and 4 for the flow formulations without and with valid inequalities, respectively, and closes the remaining five to within 0.43\%. It also returns a feasible solution on every instance, whereas the flow formulations do so on only 10 and 14 instances, at average gaps of 7.79\% and 2.83\% against 0.03\% for branch and price. 

The valid inequalities are effective, tightening the lower bound from typically 3-8\% below the best known value to within 1-2\%, and raising the number of instances solved to optimality from 1 to 4. Even so, on all but a handful of instances the bound reached after three hours remains weaker than the root node bound of the set covering formulation, which is available within seconds. Since the multi-commodity flow formulation already fails on the smallest instance class, we conclude that it is not a tractable approach to railway crew scheduling, and that branch and price is the method of choice.

\subsection{Branch and Price on Multi-Depot Instances}
\label{subsec:experiments_multi_depot}

We repeat the full branch-and-price experiments of Section~\ref{subsec:experiments_branchPrice} on the multi-depot variants of the benchmark instances described in Section~\ref{subsec:synthetic_generation}, using the same branching, primal heuristic, and reduced cost fixing settings. Table~\ref{tab:solver_comparison_multi_depot} reports the resulting bounds, gaps, and time breakdown, in the same format as Table~\ref{tab:solver_comparison}. Since every single-depot connection is also a multi-depot connection, the multi-depot instances are relaxations of their single-depot counterparts. Indeed the best known solution is at least as cheap on every multi-depot instance in our library, and strictly cheaper on all but one (the smallest instance). Optimality gaps are smaller than on the corresponding single-depot instances across all instance sizes, for example 0.52\% against 0.89\% on medium instances. This is driven primarily by the quality of the primal solutions rather than by an improved dual bound: the lower bound remains close to the root node value across all instance sizes, ranging from 100.01\% to 100.36\%. The time breakdown follows the same pattern as the single-depot results: the restricted master problem and the primal heuristics dominate computation time, while branching and reduced cost fixing contribute comparatively little. As on the single-depot instances, reduced cost fixing has no discernible effect on solution quality, but on the tiny multi-depot instances it reaches the same gap in substantially fewer nodes, 417 against 968.

\input{tables/tableSolverMultipleDepots}

The best known bounds and solutions for the multi-depot instances, obtained after additionally applying the postprocessing heuristic with the best configuration identified in Section~\ref{subsec:experiments_postprocess} (25\% of duties fixed, random selection), are reported alongside the single-depot results in Tables~\ref{tab:best_tiny}-\ref{tab:best_large} in Appendix~\ref{app:bestKnown}.

\section{Conclusion}
\label{sec:conclusion}

We conducted a systematic computational study of branch and price for railway crew scheduling, evaluating the contribution of its main components. Overall, average optimality gaps below 1\% can be obtained across all instance sizes. Primal heuristics are the main driver of solution quality: diving with backtracking produces solutions within approximately 1\% of the root node lower bound across all instance sizes, substantially outperforming a restricted master heuristic. A postprocessing heuristic reduces the upper bound by a further 0.17\% on average. Semi-strong branching offers a good balance between node quality and computation time per node, but branching contributes little to improving the dual bound, particularly on larger instances, and reduced cost fixing has negligible impact on medium and large instances. For practitioners, the overhead of a full branch-and-price scheme therefore does not appear to be warranted. Solving the root node and applying diving and postprocessing heuristics is the more practical course of action, since the additional branching effort yields only marginal improvements in the lower bound.

To support this study, we introduced a benchmark library of 100 instances for railway crew scheduling, generated by mimicking all the steps of a typical railway planning process. The instances span a wide range of sizes, are validated against historical data from Netherlands Railways, and are publicly available together with best known solutions and a solution checker. We hope this library enables a more systematic comparison of solution methods for railway crew scheduling. The instances also form a suitable basis for other problem variants, such as those involving more complex labour regulations or fairness considerations.

Our findings appear to reflect structural differences between railway crew scheduling and vehicle routing. Instance sizes are typically larger in railway crew scheduling, and crew solutions contain many more duties than vehicle routing solutions contain routes, which means that fixing a single connection rarely has a meaningful impact on the overall solution. Combined with the periodic nature of railway timetables, which provides many alternative connections between tasks, this may explain why both branching and reduced cost fixing are relatively ineffective. The computational bottleneck is the restricted master problem rather than pricing, again in contrast to vehicle routing: the acyclic pricing graph admits efficient labelling algorithms, and pricing accounts for only a small fraction of total runtime.

The central direction for future research is improving the dual bound. More effective branching rules may help, and it remains an open question whether effective valid inequalities for the set covering formulation of railway crew scheduling can be identified. Here, the problem structure implies that many of the cutting planes developed for vehicle routing do not directly carry over. For example, inequalities enforcing route elementarity have no counterpart in railway crew scheduling, as the pricing problem is already modelled on a directed acyclic graph and the resulting paths are elementary by construction. Other classical cutting planes, such as rank-one cuts and subset-row inequalities, have mostly been used for set partitioning formulations. Because pricing is relatively cheap, even non-robust cuts that require modifications to the pricing algorithm may be worthwhile, as the overhead of adapted pricing is likely outweighed by the improvement in dual bound. Moreover, it seems worthwhile to explore the use of state-expanded networks and decision diagrams as alternative solution methods for crew scheduling. Such approaches explicitly incorporate the duty constraints in the network construction, and may be better suited to leveraging the machinery of commercial integer programming solvers. Finally, it would be interesting to investigate to what extent our findings hold up in settings with more sophisticated labour regulations.

\section*{Acknowledgements}

We thank three anonymous reviewers for their suggestions that helped improve the manuscript. We also thank Jonas Brenker for critically checking earlier versions of the online repository.

\bibliographystyle{abbrvnat}
\bibliography{references}

\begin{appendices}

\section{Labelling Algorithm}
\label{app:labelling}

This appendix formalises the labelling algorithm used to solve the pricing problem of Section~\ref{subsubsec:method_pricing}. Throughout, we assume the pricing problem is solved for a given combination of crew base and starting-hour window, with the graph $G$ restricted to the source and sink arcs feasible for that crew base accordingly. We describe the label definition, the extension of forward labels, the resource feasibility checks, and the dominance rule, and sketch the bidirectional extension at the end of this appendix.

\paragraph{Notation.}
Recall that the pricing problem is a resource-constrained shortest $s,t$-path problem on the directed acyclic graph $G = (V, A)$. Every arc $a \in A$ has a cost $c_a$, start time $\underline{\tau}_{a}$, and an end time $\overline{\tau}_{a}$, so that its duration equals $\overline{\tau}_{a} - \underline{\tau}_{a}$. These times are defined such that the end time of an arc coincides with the start time of its successor on any $s,t$-path, so that the durations of the arcs on such a path sum to the length of the corresponding duty. 

Every break arc $a \in A^{B}$ additionally has a break start time $\underline{\beta}_{a}$ and a break end time $\overline{\beta}_{a}$, with $\underline{\tau}_{a} \leq \underline{\beta}_{a} \leq \overline{\beta}_{a} \leq \overline{\tau}_{a}$. For brevity, we write $\pi_{a} = \underline{\beta}_{a} - \underline{\tau}_{a}$ for the time between the start of $a$ and the start of the meal break, and $\sigma_{a} = \overline{\tau}_{a} - \overline{\beta}_{a}$ for the time between the end of the meal break and the end of $a$. Since the break itself consumes time, $\pi_{a} + \sigma_{a} \leq \overline{\tau}_{a} - \underline{\tau}_{a}$.

As observed in Section~\ref{subsubsec:method_pricing}, the reduced cost~\eqref{eq:reduced_cost} decomposes over arcs. We define the reduced cost of arc $a = (i, j) \in A$ as $\tilde{c}_{a} = c_{a} - \lambda_{j}$ if $j \in K$, and $\tilde{c}_{a} = c_{a}$ if $j = t$, so that the reduced cost of a duty equals the sum of the reduced costs of its arcs. 

\paragraph{Label definition.}
A forward label $\ell$ represents a partial duty starting at the source node $s$ and ending at some task $i \in K$. It is defined as the tuple
\begin{align*}
\ell = \left(\tilde{c}_{\ell}, \omega_{\ell}, \beta_{\ell}, b_{\ell}\right),
\end{align*}
consisting of the following components:
\begin{itemize}
\item $\tilde{c}_{\ell} \in \mathbb{R}$: the accumulated reduced cost,
\item $\omega_{\ell}$: the elapsed duty length, i.e., the time between the start of the first task of the partial duty and the end of task $i$,
\item $\beta_{\ell}$: the time without break, i.e., the time elapsed since the end of the last meal break, or since the start of the partial duty if no meal break has been taken,
\item $b_{\ell} \in \{0,1\}$: an indicator that equals one if and only if a meal break satisfying the time-before-break limit has been taken.
\end{itemize}
Together, $\omega_{\ell}$ and $\beta_{\ell}$ determine the position of the most recent meal break relative to the start of the duty. The task at which the partial duty ends need not be stored, as labels are kept in a separate bucket for every node of the graph. Similarly, the sequence of tasks visited is not stored in the label itself, but recovered from predecessor pointers once a duty of negative reduced cost is identified. The algorithm is initialised with the single label $\ell_{0} = (0, 0, 0, 0)$ in the bucket of the source node.

\paragraph{Label extension.}
Let $\ell$ be a forward label in the bucket of task $i$ and let $a = (i, j) \in A$ be an outgoing arc. Extending $\ell$ along $a$ yields a label $\ell'$ in the bucket of task $j$ with $\tilde{c}_{\ell'} = \tilde{c}_{\ell} + \tilde{c}_{a}$. The meal break resources are updated before the duty length, since the time-before-break limit is evaluated relative to the start of the duty:
\begin{align}
\left(\beta_{\ell'}, b_{\ell'}\right) =
\begin{cases}
\left(\sigma_{a}, 1\right)
& \text{if } a \in A^{B} \text{ and } \omega_{\ell} + \pi_{a} \leq L^{B}, \\
\left(\beta_{\ell} + \overline{\tau}_{a} - \underline{\tau}_{a}, b_{\ell}\right)
& \text{otherwise.}
\end{cases}
\label{eq:label_break}
\end{align}
That is, a meal break is registered on arc $a$ only if the time from the start of the partial duty to the start of the break does not exceed the maximum time without break $L^{B}$. In that case the time without break is reset to the time between the end of the break and the end of arc $a$. Otherwise, the full duration of the arc is added to the time without break. Note that a duty may contain several break arcs, but that only those starting within $L^{B}$ of the start of the duty are considered, in line with the meal break rule of Section~\ref{sec:problem}. Since $\sigma_{a} \leq \overline{\tau}_{a} - \underline{\tau}_{a}$, registering a break never increases the time without break. Finally, the duty length is updated as
\begin{align}
\omega_{\ell'} = \omega_{\ell} + \overline{\tau}_{a} - \underline{\tau}_{a}.
\label{eq:label_workload}
\end{align}
Observe that a break arc encountered later in the duty overwrites an earlier one whenever it still respects the time-before-break limit. The rule therefore positions the meal break greedily as late as that limit allows, which is optimal: a later break leaves a smaller time without break for the remainder of the duty, and can never rule out a completion that an earlier break would have permitted.

\paragraph{Resource feasibility.}
A label $\ell$ is resource feasible if its duty length does not exceed the maximum duty length $L$ and its time without break does not exceed the maximum time without break $L^{B}$, that is, if
\begin{align}
\omega_{\ell} \leq L
\quad \text{and} \quad
\beta_{\ell} \leq L^{B}.
\label{eq:label_feasible}
\end{align}
Labels violating~\eqref{eq:label_feasible} are discarded upon extension. Discarding is valid, as both resources are non-decreasing along any extension: $\omega_{\ell}$ increases monotonically by~\eqref{eq:label_workload}, and $\beta_{\ell}$ can only be reset by a break arc, which requires $\omega_{\ell} + \pi_{a} \leq L^{B}$ and hence cannot occur once $\omega_{\ell} > L^{B}$.

A label ending at the sink node $t$ corresponds to a complete duty and is feasible only if, in addition to~\eqref{eq:label_feasible}, it contains a meal break, i.e., if
\begin{align}
b_{\ell} = 1.
\label{eq:label_feasible_sink}
\end{align}

\paragraph{Dominance.}
Let $\ell$ and $\ell'$ be two forward labels in the same bucket. We say that $\ell$ \emph{dominates} $\ell'$ if
\begin{align}
\tilde{c}_{\ell} \leq \tilde{c}_{\ell'},
\quad
\omega_{\ell} \leq \omega_{\ell'},
\quad
\beta_{\ell} \leq \beta_{\ell'},
\quad \text{and} \quad
b_{\ell} \geq b_{\ell'}.
\label{eq:label_dominance}
\end{align}

Dominance is valid because both labels can be extended in exactly the same ways, and every extension preserves the four inequalities: the reduced cost and the duty length are increased by the same arc-dependent constants, while a lower elapsed duty length means that $\ell$ registers a meal break on every arc on which $\ell'$ does. Feasibility is monotone in the same direction, as a shorter duty length and a shorter time without break can never rule out a completion, and having taken a break can never be worse than not having taken one. Every feasible completion of $\ell'$ is therefore also a feasible completion of $\ell$, and yields a duty of no larger reduced cost.

\paragraph{Backward labels and concatenation.}
Backward labels represent partial duties ending at the sink node $t$, and are defined and extended in a symmetric manner, the only change being that the roles of $\pi_{a}$ and $\sigma_{a}$ are interchanged in~\eqref{eq:label_break}. Both directions are halted at a common midpoint in time, after which a forward label in the bucket of task $i$ and a backward label in the bucket of task $j$ are joined along an arc $a = (i, j) \in A$ to form a complete duty.

The reduced cost and the duty length of this duty are obtained by adding the values recorded in the two halves and the contribution of arc $a$. Only the meal break requires care, as it may be located in the forward half, in the backward half, or on the joining arc itself, so that condition~\eqref{eq:label_feasible_sink} cannot be imposed on either half separately. We therefore evaluate all three positions and retain the one leaving the smallest time without break, after which the resulting duty is checked against the maximum duty length and the maximum time without break as before.

\section{Multi-Commodity Flow Formulation}
\label{app:mcf}

We now present a multi-commodity flow formulation of the crew scheduling problem that is loosely based on the crew-indexed formulation by \citet{hoffmann2019valid}. We present it for the single-depot case, but it readily generalises to the multi-depot case. We build on the notation from Appendix~\ref{app:labelling}. Additionally, we write $A^{-}_{v}$ and $A^{+}_{v}$ to denote the set of incoming and outgoing arcs of node $v \in V$, respectively. Let $|C|$ be a valid upper bound on the number of crew members in an optimal solution. Let $C = \{1, \dots, |C|\}$ denote the set of (homogeneous) crew members, indexed by $c$. For every crew member $c \in C$ and connection $a \in A$, the binary decision variable $y_{a}^{c}$ indicates whether connection $a$ is assigned in a duty to crew member $c$. Similarly, for every $c \in C$ and break arc $a \in A^B$, binary decision variable $z_{a}^{c}$ indicates whether crew member $c$ is assigned a meal break on break arc $a$. The multi-commodity flow formulation now reads:
\begin{subequations}
\begin{align}
\min \quad & \sum_{c \in C} \sum_{a \in A} c_{a} y_{a}^{c} & \label{eq:mcf_objective} \\ 
\text{s.t.} \quad & \sum_{c \in C} \sum_{a \in A^{-}_{k}} y_{a}^{c} \geq 1 && \forall k \in K \label{eq:mcf_cover} \\
& \sum_{a \in A^{-}_{k}} y_{a}^{c} = \sum_{a \in A^{+}_{k}} y_{a}^{c} && \forall c \in C,~k \in K \label{eq:mcf_flow} \\
& \sum_{a \in A^{+}_{s}} y_{a}^{c} \leq 1 && \forall c \in C \label{eq:mcf_cardinality} \\ 
& \sum_{a \in A} (\overline{\tau}_{a} - \underline{\tau}_{a}) y_{a}^{c} \leq L \sum_{a \in A^{+}_{s}} y_{a}^{c} && \forall c \in C \label{eq:mcf_dutyLength} \\ 
& \sum_{a \in A^{B}} z_{a}^{c} = \sum_{a \in A^{+}_{s}} y_{a}^{c} && \forall c \in C \label{eq:mcf_breakSelection} \\ 
& z_{a}^{c} \leq y_{a}^{c} && \forall c \in C,~a \in A^{B} \label{eq:mcf_breakCoupling} \\ 
& \sum_{a \in A^B} \underline{\beta}_{a} z_{a}^{c} - \sum_{a \in A^{+}_{s}} \underline{\tau}_{a} y_{a}^{c}  \leq L^{B} \sum_{a \in A^{+}_{s}} y_{a}^{c} && \forall c \in C \label{eq:mcf_breakTimeBefore} \\ 
& \sum_{a \in A^{-}_{t}} \overline{\tau}_{a} y_{a}^{c} - \sum_{a \in A^B} \overline{\beta}_{a} z_{a}^{c} \leq L^{B} \sum_{a \in A^{+}_{s}} y_{a}^{c} && \forall c \in C  \label{eq:mcf_breakTimeAfter} \\ 
& y_{a}^{c} \in \{0, 1\} && \forall c \in C,~a \in A \label{eq:mcf_domainArc} \\  
& z_{a}^{c} \in \{0, 1\} && \forall c \in C,~a \in A^{B}. \label{eq:mcf_domainBreak}  
\end{align}
\label{eq:mcf}%
\end{subequations}
The objective~\eqref{eq:mcf_objective} is to minimise the sum of connection costs. Constraints~\eqref{eq:mcf_cover} specify that every task is covered at least once. Flow balance constraints~\eqref{eq:mcf_flow} and cardinality constraints~\eqref{eq:mcf_cardinality} jointly enforce that every crew member is assigned at most one proper path through the network. Constraints~\eqref{eq:mcf_dutyLength} model the duty length, while constraints~\eqref{eq:mcf_breakSelection} and~\eqref{eq:mcf_breakCoupling} require that every active crew member is assigned exactly one meal break on an active break arc. The bounds on the time before and after the break are modelled by constraints~\eqref{eq:mcf_breakTimeBefore} and~\eqref{eq:mcf_breakTimeAfter}, respectively. Constraints~\eqref{eq:mcf_domainArc}-\eqref{eq:mcf_domainBreak} model the binary domain of the decision variables.

\paragraph{Symmetry-breaking.}
Formulation~\eqref{eq:mcf} suffers from symmetry: the $|C|$ crew members are homogeneous, so every feasible solution can be represented in $|C|!$ different ways. We break part of this symmetry in two ways.

First, following \citet{hoffmann2019valid}, we require crew member $c$ to be active whenever crew member $c+1$ is active, by imposing the ordering constraints
\begin{align}
    \sum_{a \in A_{s}^{+}} y_{a}^{c} \geq \sum_{a \in A_{s}^{+}} y_{a}^{c + 1} \quad \forall c \in C \setminus \{ |C| \}. \label{eq:mcf_ordering}
\end{align}

Second, we assign a set of pairwise incompatible tasks to a fixed set of crew members. Let $Q \subseteq K$ be a set of tasks no two of which can be covered by the same crew member in a feasible duty. Ordering these tasks arbitrarily as $Q = \{ q_{1}, \dots, q_{|Q|} \}$, we assign them to the first $|Q|$ crew members through the constraints
\begin{align}
    \sum_{a \in A^{-}_{q_c}} y^{c}_{a} = 1 \quad \forall c = 1, \dots, |Q|. \label{eq:mcf_clique_task}
\end{align}
This assignment is consistent with~\eqref{eq:mcf_ordering}, and its strength grows with the size of $Q$. We construct a large set $Q$ as follows. Call a task \emph{active} at time $\tau$ if it starts at or before $\tau$ and ends after $\tau$. For every task start time $\tau$, let $Q_{\tau}$ consist of all tasks that are active at $\tau$ together with all tasks that are active at $\tau + L + 1$. Any two tasks in $Q_{\tau}$ are indeed incompatible: two tasks that are active at the same point in time overlap, whereas a duty covering a task active at $\tau$ and a task active at $\tau + L + 1$ would have a length exceeding $L$. Since the set of active tasks only changes at task start times, it suffices to enumerate these, and we set $Q$ equal to the largest set $Q_{\tau}$ obtained this way. Finally, since crew member $c \leq |Q|$ covers task $q_{c}$, we fix $y^{c}_{a} = 0$ for every arc $a$ that cannot appear in a feasible duty covering $q_{c}$, i.e., every arc whose start or end time violates the maximum duty length when combined with $q_{c}$.

\paragraph{Valid inequalities.}
Constraints~\eqref{eq:mcf_dutyLength}, \eqref{eq:mcf_breakTimeBefore}, and~\eqref{eq:mcf_breakTimeAfter} enforce the duty length and meal break rules through a single aggregate inequality per crew member, which is weak in the linear programming relaxation. We strengthen the formulation with a family of inequalities that impose these rules at every individual point in time, loosely inspired by the clique-based formulation of \citet{breugem2022column}.

Consider a crew member $c \in C$ and a point in time $\tau$. If the duty of $c$ starts at or before $\tau$, it must end at or before $\tau + L$, giving the duty length inequalities
\begin{align}
\sum_{a \in A_{s}^{+} : \underline{\tau}_{a} \leq \tau} y_{a}^{c} + \sum_{a \in A_{t}^{-} : \overline{\tau}_{a} > \tau + L} y_{a}^{c} \leq \sum_{a \in A_{s}^{+}} y_{a}^{c}. \label{eq:mcf_vi_length}
\end{align}
The same argument applies to the meal break. If the duty of $c$ starts at or before $\tau$, its break must start at or before $\tau + L^{B}$, and if its break ends at or before $\tau$, the duty must end at or before $\tau + L^{B}$. This yields
\begin{align}
\sum_{a \in A_{s}^{+} : \underline{\tau}_{a} \leq \tau} y_{a}^{c} + \sum_{a \in A^{B} : \underline{\beta}_{a} > \tau + L^{B}} z_{a}^{c} &\leq \sum_{a \in A_{s}^{+}} y_{a}^{c}, \label{eq:mcf_vi_breakBefore} \\
\sum_{a \in A^{B} : \overline{\beta}_{a} \leq \tau} z_{a}^{c} + \sum_{a \in A_{t}^{-} : \overline{\tau}_{a} > \tau + L^{B}} y_{a}^{c} &\leq \sum_{a \in A_{s}^{+}} y_{a}^{c}. \label{eq:mcf_vi_breakAfter}
\end{align}
If crew member $c$ is not active, the right-hand side is zero and all terms on the left-hand side vanish by flow conservation. If $c$ is active, exactly one source arc, one sink arc, and one break arc are selected by~\eqref{eq:mcf_breakSelection}, so that each sum on the left-hand side is at most one, and the inequality forbids precisely the combinations of events that violate the meal break rule. We include the above valid inequalities for all unique task start times and all crew members.

\section{Best Known Bounds}
\label{app:bestKnown}

\renewcommand{\thetable}{C\arabic{table}}
\setcounter{table}{0}

Tables~\ref{tab:best_tiny}-\ref{tab:best_large} report the root node lower bound, best known lower bound, best known upper bound, optimality gap, and number of duties in the best known solutions for the tiny, small, medium, and large instances, respectively. The best known lower bound for an instance is the maximum final lower bound observed across the branching-rule experiments of Section~\ref{subsec:experiments_branching} and the full branch-and-price runs of Section~\ref{subsec:experiments_branchPrice}. The best known upper bound is the minimum cost of any feasible solution found across the primal heuristic experiments of Section~\ref{subsec:experiments_primal}, the full branch-and-price runs of Section~\ref{subsec:experiments_branchPrice}, and the postprocessing heuristic of Section~\ref{subsec:experiments_postprocess}. The corresponding solution files are available in the benchmark repository at \url{https://github.com/BartvanRossum/rcsp-benchmark}.

\input{tables/tableBestKnownTiny}

\input{tables/tableBestKnownSmall}

\input{tables/tableBestKnownMedium}

\input{tables/tableBestKnownLarge}

\section{Sensitivity Analysis of Column Generation Parameters}
\label{sec:parameter}

\renewcommand{\thefigure}{D\arabic{figure}}
\setcounter{figure}{0}

We conduct a parameter sensitivity analysis to assess the effect of the column selection, column management, and dual smoothing strategies presented in Section~\ref{subsubsec:method_accelerate}. We consider three possible parameter settings for each strategy. For column selection, we consider adding columns that are at least 25, 50, or 75\% task-disjoint. For column management, we consider removing columns that have been inactive in the last 10, 20, or 50 iterations. Finally, for dual smoothing, we consider (i) no smoothing ($\alpha=0$), (ii) dual smoothing with a fixed smoothing parameter $\alpha=0.5$, and (iii) dual smoothing with a self-tuning smoothing parameter $\alpha$. For all $3 \times 3 \times 3 = 27$ possible configurations, we solve the root node for all 25 medium-sized instances and report the average solution time (in seconds) and the average number of column generation iterations.

\input{figs/parameterTuning}

Each marker in Figure~\ref{fig:parameter_tuning} represents a single parameter configuration, averaged over the 25 medium-sized instances, and markers are coloured and shaped within each panel to distinguish the three settings of the parameter being varied. Dual smoothing has the clearest effect: solving without smoothing ($\alpha = 0$) is noticeably slower than either smoothed variant, despite a similar number of iterations, while a fixed and a self-tuning smoothing parameter perform comparably. The column selection threshold has a more moderate effect, with a looser threshold of 25\% requiring more iterations and time than the two stricter thresholds, whereas the column removal threshold has little visible impact on either measure. Overall, the parameter setting combining a 50\% task-disjointness threshold, column removal after 10 iterations, and a self-tuning dual smoothing parameter proves to be the fastest on average, and is the setting we use in our computational experiments.

\end{appendices}

\end{document}

%% file: figs/generation.tex
\begin{figure}[htbp]
\centering
\begin{tikzpicture}[
    node distance=0.4cm,
    box/.style={
        rectangle, rounded corners=2pt,
        draw=black, fill=gray!15,
        minimum width=1.75cm, align=center,
        font=\small, inner sep=4pt, minimum height=1.0cm
    },
    arr/.style={->, >=stealth, thick}
]
\node[box] (infra)   {Infrastructure\\network};
\node[box, right=of infra]  (demand)  {Passenger\\demand};
\node[box, right=of demand] (line)    {Line plan};
\node[box, right=of line]   (tt)      {Timetable};
\node[box, right=of tt]     (rs)      {Rolling stock\\schedule};
\node[box, right=of rs]     (tasks)   {Tasks \&\\connections};

\draw[arr] (infra)  -- (demand);
\draw[arr] (demand) -- (line);
\draw[arr] (line)   -- (tt);
\draw[arr] (tt)     -- (rs);
\draw[arr] (rs)     -- (tasks);
\end{tikzpicture}
\caption{Instance generation pipeline.}
\label{fig:generation}
\end{figure}

%% file: tables/instanceTable.tex
\begin{table}[htbp]
\centering
\caption{Benchmark instance statistics per group.}
\label{tab:instances}
\begin{adjustbox}{max width=\textwidth}
\begin{tabular}{lrrrrrr}
\toprule
& & & & \multicolumn{2}{c}{Connections} & \\
\cmidrule(lr){5-6}
Group & Instances & Stations & Tasks & Single-depot & Multi-depot & Task duration (min) \\
\midrule
Tiny & 25 & 3.7 & 276 & 2,843 & 3,205 & 33.9 \\
Small & 25 & 7.6 & 715 & 9,745 & 10,857 & 30.1 \\
Medium & 25 & 12.6 & 1,414 & 24,357 & 26,744 & 29.7 \\
Large & 25 & 17.8 & 2,440 & 53,339 & 57,612 & 29.8 \\
\bottomrule
\end{tabular}
\end{adjustbox}
\end{table}

%% file: figs/instancePlot.tex
\begin{figure}[htbp]
\centering
\begin{subfigure}[b]{0.45\textwidth}
\centering
\begin{tikzpicture}[x=1cm, y=1cm]
  \draw[black, line width=0.8pt] (4.717,1.099) -- (5.542,1.109);
  \draw[black, line width=0.8pt] (4.717,1.099) -- (3.770,0.000);
  \draw[black, line width=0.8pt] (2.009,3.673) -- (2.577,6.000);
  \draw[black, line width=0.8pt] (5.542,1.109) -- (3.770,0.000);
  \draw[black, line width=0.8pt] (3.770,0.000) -- (3.453,0.804);
  \draw[black, line width=0.8pt] (4.717,1.099) -- (3.060,3.291);
  \draw[black, line width=0.8pt] (2.009,3.673) -- (3.060,3.291);
  \draw[black, line width=0.8pt] (3.770,0.000) -- (2.734,0.254);
  \draw[black, line width=0.8pt] (2.009,3.673) -- (1.423,2.889);
  \draw[black, line width=0.8pt] (0.000,1.914) -- (1.423,2.889);
  \draw[black, line width=0.8pt] (0.000,1.914) -- (0.198,0.945);
  \draw[black, line width=0.8pt] (2.734,0.254) -- (0.198,0.945);
  \draw[black, line width=0.8pt] (2.009,3.673) -- (3.553,1.589);
  \draw[black, line width=0.8pt] (3.453,0.804) -- (3.553,1.589);
  \node[draw=black, fill=white, circle, minimum size=4.5pt, inner sep=0pt, line width=0.4pt] at (4.717,1.099) {};
  \node[draw=black, fill=white, circle, minimum size=4.5pt, inner sep=0pt, line width=0.4pt] at (5.542,1.109) {};
  \node[draw=black, fill=white, circle, minimum size=4.5pt, inner sep=0pt, line width=0.4pt] at (2.009,3.673) {};
  \node[draw=black, fill=white, circle, minimum size=4.5pt, inner sep=0pt, line width=0.4pt] at (2.577,6.000) {};
  \node[draw=black, fill=white, circle, minimum size=4.5pt, inner sep=0pt, line width=0.4pt] at (0.000,1.914) {};
  \node[draw=black, fill=gray!75, circle, minimum size=3pt, inner sep=0pt, line width=0.4pt] at (3.453,0.804) {};
  \node[draw=black, fill=gray!75, circle, minimum size=3pt, inner sep=0pt, line width=0.4pt] at (3.060,3.291) {};
  \node[draw=black, fill=gray!75, circle, minimum size=3pt, inner sep=0pt, line width=0.4pt] at (2.734,0.254) {};
  \node[draw=black, fill=gray!75, circle, minimum size=3pt, inner sep=0pt, line width=0.4pt] at (1.423,2.889) {};
  \node[draw=black, fill=gray!75, circle, minimum size=3pt, inner sep=0pt, line width=0.4pt] at (0.198,0.945) {};
  \node[draw=black, fill=gray!75, circle, minimum size=3pt, inner sep=0pt, line width=0.4pt] at (3.553,1.589) {};
  \node[draw=black, fill=white, rectangle, minimum size=6pt, inner sep=0pt, line width=0.6pt] at (3.770,0.000) {};
\end{tikzpicture}
\caption{\texttt{s6\_t482}}
\end{subfigure}
\hspace{1cm}
\begin{subfigure}[b]{0.45\textwidth}
\centering
\begin{tikzpicture}[x=1cm, y=1cm]
  \draw[black, line width=0.8pt] (2.017,2.495) -- (2.598,3.261);
  \draw[black, line width=0.8pt] (1.418,3.554) -- (2.598,3.261);
  \draw[black, line width=0.8pt] (2.009,1.796) -- (3.455,2.016);
  \draw[black, line width=0.8pt] (3.893,2.531) -- (2.598,3.261);
  \draw[black, line width=0.8pt] (2.210,0.000) -- (3.455,2.016);
  \draw[black, line width=0.8pt] (3.893,2.531) -- (3.849,2.152);
  \draw[black, line width=0.8pt] (3.455,2.016) -- (3.849,2.152);
  \draw[black, line width=0.8pt] (2.017,2.495) -- (1.499,2.779);
  \draw[black, line width=0.8pt] (1.418,3.554) -- (1.499,2.779);
  \draw[black, line width=0.8pt] (2.009,1.796) -- (1.884,2.201);
  \draw[black, line width=0.8pt] (2.017,2.495) -- (1.884,2.201);
  \draw[black, line width=0.8pt] (2.210,0.000) -- (2.129,0.294);
  \draw[black, line width=0.8pt] (2.009,1.796) -- (1.852,1.574);
  \draw[black, line width=0.8pt] (2.129,0.294) -- (1.852,1.574);
  \draw[black, line width=0.8pt] (2.009,1.796) -- (1.447,2.177);
  \draw[black, line width=0.8pt] (0.000,2.868) -- (1.447,2.177);
  \draw[black, line width=0.8pt] (1.418,3.554) -- (0.442,3.336);
  \draw[black, line width=0.8pt] (0.000,2.868) -- (0.345,3.079);
  \draw[black, line width=0.8pt] (0.442,3.336) -- (0.345,3.079);
  \draw[black, line width=0.8pt] (3.893,2.531) -- (5.856,3.176);
  \draw[black, line width=0.8pt] (6.000,3.579) -- (5.856,3.176);
  \node[draw=black, fill=white, circle, minimum size=4.5pt, inner sep=0pt, line width=0.4pt] at (3.893,2.531) {};
  \node[draw=black, fill=white, circle, minimum size=4.5pt, inner sep=0pt, line width=0.4pt] at (6.000,3.579) {};
  \node[draw=black, fill=white, circle, minimum size=4.5pt, inner sep=0pt, line width=0.4pt] at (2.017,2.495) {};
  \node[draw=black, fill=white, circle, minimum size=4.5pt, inner sep=0pt, line width=0.4pt] at (1.418,3.554) {};
  \node[draw=black, fill=white, circle, minimum size=4.5pt, inner sep=0pt, line width=0.4pt] at (2.210,0.000) {};
  \node[draw=black, fill=white, circle, minimum size=4.5pt, inner sep=0pt, line width=0.4pt] at (2.598,3.261) {};
  \node[draw=black, fill=white, circle, minimum size=4.5pt, inner sep=0pt, line width=0.4pt] at (3.455,2.016) {};
  \node[draw=black, fill=white, circle, minimum size=4.5pt, inner sep=0pt, line width=0.4pt] at (0.000,2.868) {};
  \node[draw=black, fill=gray!75, circle, minimum size=3pt, inner sep=0pt, line width=0.4pt] at (3.849,2.152) {};
  \node[draw=black, fill=gray!75, circle, minimum size=3pt, inner sep=0pt, line width=0.4pt] at (1.499,2.779) {};
  \node[draw=black, fill=gray!75, circle, minimum size=3pt, inner sep=0pt, line width=0.4pt] at (1.884,2.201) {};
  \node[draw=black, fill=gray!75, circle, minimum size=3pt, inner sep=0pt, line width=0.4pt] at (2.129,0.294) {};
  \node[draw=black, fill=gray!75, circle, minimum size=3pt, inner sep=0pt, line width=0.4pt] at (1.852,1.574) {};
  \node[draw=black, fill=gray!75, circle, minimum size=3pt, inner sep=0pt, line width=0.4pt] at (1.447,2.177) {};
  \node[draw=black, fill=gray!75, circle, minimum size=3pt, inner sep=0pt, line width=0.4pt] at (0.442,3.336) {};
  \node[draw=black, fill=gray!75, circle, minimum size=3pt, inner sep=0pt, line width=0.4pt] at (0.345,3.079) {};
  \node[draw=black, fill=gray!75, circle, minimum size=3pt, inner sep=0pt, line width=0.4pt] at (5.856,3.176) {};
  \node[draw=black, fill=white, rectangle, minimum size=6pt, inner sep=0pt, line width=0.6pt] at (2.009,1.796) {};
\end{tikzpicture}
\caption{\texttt{s9\_t991}}
\end{subfigure}
\\[1em]
\begin{subfigure}[b]{0.45\textwidth}
\centering
\begin{tikzpicture}[x=1cm, y=1cm]
  \draw[black, line width=0.8pt] (0.000,0.390) -- (0.927,0.000);
  \draw[black, line width=0.8pt] (0.927,0.000) -- (1.436,0.838);
  \draw[black, line width=0.8pt] (2.121,1.316) -- (1.436,0.838);
  \draw[black, line width=0.8pt] (2.121,1.316) -- (1.990,1.875);
  \draw[black, line width=0.8pt] (1.990,1.875) -- (2.011,2.466);
  \draw[black, line width=0.8pt] (4.106,2.750) -- (4.621,2.258);
  \draw[black, line width=0.8pt] (5.250,5.007) -- (4.192,4.462);
  \draw[black, line width=0.8pt] (2.121,1.316) -- (4.051,2.030);
  \draw[black, line width=0.8pt] (4.621,2.258) -- (4.051,2.030);
  \draw[black, line width=0.8pt] (4.106,2.750) -- (3.522,2.827);
  \draw[black, line width=0.8pt] (0.000,0.390) -- (1.694,1.087);
  \draw[black, line width=0.8pt] (2.121,1.316) -- (1.694,1.087);
  \draw[black, line width=0.8pt] (2.011,2.466) -- (3.056,2.852);
  \draw[black, line width=0.8pt] (3.522,2.827) -- (3.056,2.852);
  \draw[black, line width=0.8pt] (2.011,2.466) -- (1.954,3.543);
  \draw[black, line width=0.8pt] (1.754,4.184) -- (1.954,3.543);
  \draw[black, line width=0.8pt] (4.106,2.750) -- (4.988,4.485);
  \draw[black, line width=0.8pt] (5.250,5.007) -- (4.988,4.485);
  \draw[black, line width=0.8pt] (1.754,4.184) -- (1.208,5.335);
  \draw[black, line width=0.8pt] (0.836,6.000) -- (1.208,5.335);
  \draw[black, line width=0.8pt] (4.106,2.750) -- (4.071,3.091);
  \draw[black, line width=0.8pt] (4.192,4.462) -- (4.071,3.091);
  \draw[black, line width=0.8pt] (4.106,2.750) -- (2.421,3.428);
  \draw[black, line width=0.8pt] (1.754,4.184) -- (2.421,3.428);
  \draw[black, line width=0.8pt] (4.192,4.462) -- (3.964,4.691);
  \draw[black, line width=0.8pt] (2.121,1.316) -- (3.049,2.138);
  \draw[black, line width=0.8pt] (4.106,2.750) -- (3.049,2.138);
  \draw[black, line width=0.8pt] (1.754,4.184) -- (2.395,4.356);
  \draw[black, line width=0.8pt] (3.964,4.691) -- (2.395,4.356);
  \node[draw=black, fill=white, circle, minimum size=4.5pt, inner sep=0pt, line width=0.4pt] at (0.000,0.390) {};
  \node[draw=black, fill=white, circle, minimum size=4.5pt, inner sep=0pt, line width=0.4pt] at (2.121,1.316) {};
  \node[draw=black, fill=white, circle, minimum size=4.5pt, inner sep=0pt, line width=0.4pt] at (1.990,1.875) {};
  \node[draw=black, fill=white, circle, minimum size=4.5pt, inner sep=0pt, line width=0.4pt] at (0.927,0.000) {};
  \node[draw=black, fill=white, circle, minimum size=4.5pt, inner sep=0pt, line width=0.4pt] at (4.621,2.258) {};
  \node[draw=black, fill=white, circle, minimum size=4.5pt, inner sep=0pt, line width=0.4pt] at (1.754,4.184) {};
  \node[draw=black, fill=white, circle, minimum size=4.5pt, inner sep=0pt, line width=0.4pt] at (5.250,5.007) {};
  \node[draw=black, fill=white, circle, minimum size=4.5pt, inner sep=0pt, line width=0.4pt] at (4.192,4.462) {};
  \node[draw=black, fill=white, circle, minimum size=4.5pt, inner sep=0pt, line width=0.4pt] at (2.011,2.466) {};
  \node[draw=black, fill=white, circle, minimum size=4.5pt, inner sep=0pt, line width=0.4pt] at (1.436,0.838) {};
  \node[draw=black, fill=white, circle, minimum size=4.5pt, inner sep=0pt, line width=0.4pt] at (0.836,6.000) {};
  \node[draw=black, fill=gray!75, circle, minimum size=3pt, inner sep=0pt, line width=0.4pt] at (4.051,2.030) {};
  \node[draw=black, fill=gray!75, circle, minimum size=3pt, inner sep=0pt, line width=0.4pt] at (3.522,2.827) {};
  \node[draw=black, fill=gray!75, circle, minimum size=3pt, inner sep=0pt, line width=0.4pt] at (1.694,1.087) {};
  \node[draw=black, fill=gray!75, circle, minimum size=3pt, inner sep=0pt, line width=0.4pt] at (3.056,2.852) {};
  \node[draw=black, fill=gray!75, circle, minimum size=3pt, inner sep=0pt, line width=0.4pt] at (1.954,3.543) {};
  \node[draw=black, fill=gray!75, circle, minimum size=3pt, inner sep=0pt, line width=0.4pt] at (4.988,4.485) {};
  \node[draw=black, fill=gray!75, circle, minimum size=3pt, inner sep=0pt, line width=0.4pt] at (1.208,5.335) {};
  \node[draw=black, fill=gray!75, circle, minimum size=3pt, inner sep=0pt, line width=0.4pt] at (4.071,3.091) {};
  \node[draw=black, fill=gray!75, circle, minimum size=3pt, inner sep=0pt, line width=0.4pt] at (2.421,3.428) {};
  \node[draw=black, fill=gray!75, circle, minimum size=3pt, inner sep=0pt, line width=0.4pt] at (3.964,4.691) {};
  \node[draw=black, fill=gray!75, circle, minimum size=3pt, inner sep=0pt, line width=0.4pt] at (3.049,2.138) {};
  \node[draw=black, fill=gray!75, circle, minimum size=3pt, inner sep=0pt, line width=0.4pt] at (2.395,4.356) {};
  \node[draw=black, fill=white, rectangle, minimum size=6pt, inner sep=0pt, line width=0.6pt] at (4.106,2.750) {};
\end{tikzpicture}
\caption{\texttt{s12\_t1459}}
\end{subfigure}
\hspace{1cm}
\begin{subfigure}[b]{0.45\textwidth}
\centering
\begin{tikzpicture}[x=1cm, y=1cm]
  \draw[black, line width=0.8pt] (2.967,0.496) -- (2.466,0.396);
  \draw[black, line width=0.8pt] (2.373,0.809) -- (2.466,0.396);
  \draw[black, line width=0.8pt] (2.657,1.724) -- (1.781,1.504);
  \draw[black, line width=0.8pt] (1.212,0.000) -- (0.841,0.446);
  \draw[black, line width=0.8pt] (3.632,0.377) -- (4.963,0.889);
  \draw[black, line width=0.8pt] (4.963,0.889) -- (5.233,1.344);
  \draw[black, line width=0.8pt] (4.940,1.885) -- (5.233,1.344);
  \draw[black, line width=0.8pt] (0.841,0.446) -- (1.781,1.504);
  \draw[black, line width=0.8pt] (2.967,0.496) -- (2.373,0.809);
  \draw[black, line width=0.8pt] (2.657,1.724) -- (1.944,2.485);
  \draw[black, line width=0.8pt] (5.825,2.602) -- (6.000,2.179);
  \draw[black, line width=0.8pt] (2.466,0.396) -- (2.069,0.229);
  \draw[black, line width=0.8pt] (2.657,1.724) -- (2.416,1.260);
  \draw[black, line width=0.8pt] (2.967,0.496) -- (2.416,1.260);
  \draw[black, line width=0.8pt] (2.967,0.496) -- (3.428,0.472);
  \draw[black, line width=0.8pt] (3.632,0.377) -- (3.428,0.472);
  \draw[black, line width=0.8pt] (3.632,0.377) -- (4.587,0.889);
  \draw[black, line width=0.8pt] (5.233,1.344) -- (4.587,0.889);
  \draw[black, line width=0.8pt] (0.000,0.780) -- (0.184,0.757);
  \draw[black, line width=0.8pt] (0.841,0.446) -- (0.184,0.757);
  \draw[black, line width=0.8pt] (4.940,1.885) -- (5.027,2.301);
  \draw[black, line width=0.8pt] (5.825,2.602) -- (5.027,2.301);
  \draw[black, line width=0.8pt] (0.841,0.446) -- (1.199,0.491);
  \draw[black, line width=0.8pt] (2.373,0.809) -- (2.212,1.273);
  \draw[black, line width=0.8pt] (1.781,1.504) -- (2.212,1.273);
  \draw[black, line width=0.8pt] (5.233,1.344) -- (5.690,2.173);
  \draw[black, line width=0.8pt] (6.000,2.179) -- (5.690,2.173);
  \draw[black, line width=0.8pt] (2.466,0.396) -- (2.047,0.509);
  \draw[black, line width=0.8pt] (1.199,0.491) -- (1.570,0.452);
  \draw[black, line width=0.8pt] (2.047,0.509) -- (1.570,0.452);
  \draw[black, line width=0.8pt] (1.944,2.485) -- (1.598,1.801);
  \draw[black, line width=0.8pt] (1.781,1.504) -- (1.598,1.801);
  \draw[black, line width=0.8pt] (1.212,0.000) -- (2.035,0.062);
  \draw[black, line width=0.8pt] (2.069,0.229) -- (2.035,0.062);
  \node[draw=black, fill=white, circle, minimum size=4.5pt, inner sep=0pt, line width=0.4pt] at (2.373,0.809) {};
  \node[draw=black, fill=white, circle, minimum size=4.5pt, inner sep=0pt, line width=0.4pt] at (3.632,0.377) {};
  \node[draw=black, fill=white, circle, minimum size=4.5pt, inner sep=0pt, line width=0.4pt] at (2.466,0.396) {};
  \node[draw=black, fill=white, circle, minimum size=4.5pt, inner sep=0pt, line width=0.4pt] at (1.212,0.000) {};
  \node[draw=black, fill=white, circle, minimum size=4.5pt, inner sep=0pt, line width=0.4pt] at (2.657,1.724) {};
  \node[draw=black, fill=white, circle, minimum size=4.5pt, inner sep=0pt, line width=0.4pt] at (0.841,0.446) {};
  \node[draw=black, fill=white, circle, minimum size=4.5pt, inner sep=0pt, line width=0.4pt] at (1.944,2.485) {};
  \node[draw=black, fill=white, circle, minimum size=4.5pt, inner sep=0pt, line width=0.4pt] at (0.000,0.780) {};
  \node[draw=black, fill=white, circle, minimum size=4.5pt, inner sep=0pt, line width=0.4pt] at (1.781,1.504) {};
  \node[draw=black, fill=white, circle, minimum size=4.5pt, inner sep=0pt, line width=0.4pt] at (4.963,0.889) {};
  \node[draw=black, fill=white, circle, minimum size=4.5pt, inner sep=0pt, line width=0.4pt] at (4.940,1.885) {};
  \node[draw=black, fill=white, circle, minimum size=4.5pt, inner sep=0pt, line width=0.4pt] at (5.233,1.344) {};
  \node[draw=black, fill=white, circle, minimum size=4.5pt, inner sep=0pt, line width=0.4pt] at (5.825,2.602) {};
  \node[draw=black, fill=gray!75, circle, minimum size=3pt, inner sep=0pt, line width=0.4pt] at (6.000,2.179) {};
  \node[draw=black, fill=gray!75, circle, minimum size=3pt, inner sep=0pt, line width=0.4pt] at (2.069,0.229) {};
  \node[draw=black, fill=gray!75, circle, minimum size=3pt, inner sep=0pt, line width=0.4pt] at (2.416,1.260) {};
  \node[draw=black, fill=gray!75, circle, minimum size=3pt, inner sep=0pt, line width=0.4pt] at (3.428,0.472) {};
  \node[draw=black, fill=gray!75, circle, minimum size=3pt, inner sep=0pt, line width=0.4pt] at (4.587,0.889) {};
  \node[draw=black, fill=gray!75, circle, minimum size=3pt, inner sep=0pt, line width=0.4pt] at (0.184,0.757) {};
  \node[draw=black, fill=gray!75, circle, minimum size=3pt, inner sep=0pt, line width=0.4pt] at (5.027,2.301) {};
  \node[draw=black, fill=gray!75, circle, minimum size=3pt, inner sep=0pt, line width=0.4pt] at (1.199,0.491) {};
  \node[draw=black, fill=gray!75, circle, minimum size=3pt, inner sep=0pt, line width=0.4pt] at (2.212,1.273) {};
  \node[draw=black, fill=gray!75, circle, minimum size=3pt, inner sep=0pt, line width=0.4pt] at (5.690,2.173) {};
  \node[draw=black, fill=gray!75, circle, minimum size=3pt, inner sep=0pt, line width=0.4pt] at (2.047,0.509) {};
  \node[draw=black, fill=gray!75, circle, minimum size=3pt, inner sep=0pt, line width=0.4pt] at (1.570,0.452) {};
  \node[draw=black, fill=gray!75, circle, minimum size=3pt, inner sep=0pt, line width=0.4pt] at (1.598,1.801) {};
  \node[draw=black, fill=gray!75, circle, minimum size=3pt, inner sep=0pt, line width=0.4pt] at (2.035,0.062) {};
  \node[draw=black, fill=white, rectangle, minimum size=6pt, inner sep=0pt, line width=0.6pt] at (2.967,0.496) {};
\end{tikzpicture}
\caption{\texttt{s14\_t1755}}
\end{subfigure}
\\[1em]
\begin{subfigure}[b]{0.45\textwidth}
\centering
\begin{tikzpicture}[x=1cm, y=1cm]
  \draw[black, line width=0.8pt] (3.803,3.485) -- (3.589,2.759);
  \draw[black, line width=0.8pt] (3.803,3.485) -- (3.269,4.404);
  \draw[black, line width=0.8pt] (3.269,4.404) -- (1.897,4.409);
  \draw[black, line width=0.8pt] (1.897,4.409) -- (1.341,4.844);
  \draw[black, line width=0.8pt] (1.341,4.844) -- (0.570,4.474);
  \draw[black, line width=0.8pt] (1.534,3.394) -- (1.124,2.939);
  \draw[black, line width=0.8pt] (0.570,4.474) -- (0.000,5.585);
  \draw[black, line width=0.8pt] (2.209,1.409) -- (2.689,1.291);
  \draw[black, line width=0.8pt] (2.396,0.697) -- (2.689,1.291);
  \draw[black, line width=0.8pt] (2.040,0.388) -- (2.396,0.697);
  \draw[black, line width=0.8pt] (3.809,0.898) -- (4.511,0.000);
  \draw[black, line width=0.8pt] (4.511,0.000) -- (6.000,0.632);
  \draw[black, line width=0.8pt] (3.809,0.898) -- (6.000,0.632);
  \draw[black, line width=0.8pt] (3.809,0.898) -- (3.946,1.500);
  \draw[black, line width=0.8pt] (3.589,2.759) -- (2.242,4.237);
  \draw[black, line width=0.8pt] (1.897,4.409) -- (2.242,4.237);
  \draw[black, line width=0.8pt] (3.809,0.898) -- (3.739,0.642);
  \draw[black, line width=0.8pt] (1.897,4.409) -- (1.801,3.803);
  \draw[black, line width=0.8pt] (1.534,3.394) -- (1.801,3.803);
  \draw[black, line width=0.8pt] (2.209,1.409) -- (1.179,2.645);
  \draw[black, line width=0.8pt] (1.124,2.939) -- (1.179,2.645);
  \draw[black, line width=0.8pt] (2.689,1.291) -- (2.874,1.837);
  \draw[black, line width=0.8pt] (3.589,2.759) -- (2.874,1.837);
  \draw[black, line width=0.8pt] (3.589,2.759) -- (3.099,2.625);
  \draw[black, line width=0.8pt] (3.589,2.759) -- (3.485,1.888);
  \draw[black, line width=0.8pt] (3.946,1.500) -- (3.485,1.888);
  \draw[black, line width=0.8pt] (2.689,1.291) -- (2.582,2.336);
  \draw[black, line width=0.8pt] (2.582,2.336) -- (2.199,2.598);
  \draw[black, line width=0.8pt] (1.534,3.394) -- (2.328,3.024);
  \draw[black, line width=0.8pt] (3.099,2.625) -- (2.328,3.024);
  \draw[black, line width=0.8pt] (2.689,1.291) -- (3.367,0.761);
  \draw[black, line width=0.8pt] (3.809,0.898) -- (3.367,0.761);
  \draw[black, line width=0.8pt] (2.926,0.055) -- (3.783,0.413);
  \draw[black, line width=0.8pt] (3.739,0.642) -- (3.783,0.413);
  \draw[black, line width=0.8pt] (0.570,4.474) -- (0.551,4.011);
  \draw[black, line width=0.8pt] (1.534,3.394) -- (1.506,2.919);
  \draw[black, line width=0.8pt] (2.199,2.598) -- (1.506,2.919);
  \draw[black, line width=0.8pt] (2.396,0.697) -- (2.624,0.553);
  \draw[black, line width=0.8pt] (2.926,0.055) -- (2.624,0.553);
  \draw[black, line width=0.8pt] (0.551,4.011) -- (1.362,3.562);
  \draw[black, line width=0.8pt] (1.534,3.394) -- (1.362,3.562);
  \node[draw=black, fill=white, circle, minimum size=4.5pt, inner sep=0pt, line width=0.4pt] at (3.803,3.485) {};
  \node[draw=black, fill=white, circle, minimum size=4.5pt, inner sep=0pt, line width=0.4pt] at (2.209,1.409) {};
  \node[draw=black, fill=white, circle, minimum size=4.5pt, inner sep=0pt, line width=0.4pt] at (3.269,4.404) {};
  \node[draw=black, fill=white, circle, minimum size=4.5pt, inner sep=0pt, line width=0.4pt] at (2.040,0.388) {};
  \node[draw=black, fill=white, circle, minimum size=4.5pt, inner sep=0pt, line width=0.4pt] at (2.396,0.697) {};
  \node[draw=black, fill=white, circle, minimum size=4.5pt, inner sep=0pt, line width=0.4pt] at (1.897,4.409) {};
  \node[draw=black, fill=white, circle, minimum size=4.5pt, inner sep=0pt, line width=0.4pt] at (2.689,1.291) {};
  \node[draw=black, fill=white, circle, minimum size=4.5pt, inner sep=0pt, line width=0.4pt] at (1.341,4.844) {};
  \node[draw=black, fill=white, circle, minimum size=4.5pt, inner sep=0pt, line width=0.4pt] at (2.926,0.055) {};
  \node[draw=black, fill=white, circle, minimum size=4.5pt, inner sep=0pt, line width=0.4pt] at (0.570,4.474) {};
  \node[draw=black, fill=white, circle, minimum size=4.5pt, inner sep=0pt, line width=0.4pt] at (3.809,0.898) {};
  \node[draw=black, fill=white, circle, minimum size=4.5pt, inner sep=0pt, line width=0.4pt] at (4.511,0.000) {};
  \node[draw=black, fill=white, circle, minimum size=4.5pt, inner sep=0pt, line width=0.4pt] at (1.534,3.394) {};
  \node[draw=black, fill=white, circle, minimum size=4.5pt, inner sep=0pt, line width=0.4pt] at (1.124,2.939) {};
  \node[draw=black, fill=white, circle, minimum size=4.5pt, inner sep=0pt, line width=0.4pt] at (6.000,0.632) {};
  \node[draw=black, fill=white, circle, minimum size=4.5pt, inner sep=0pt, line width=0.4pt] at (0.000,5.585) {};
  \node[draw=black, fill=gray!75, circle, minimum size=3pt, inner sep=0pt, line width=0.4pt] at (3.946,1.500) {};
  \node[draw=black, fill=gray!75, circle, minimum size=3pt, inner sep=0pt, line width=0.4pt] at (2.242,4.237) {};
  \node[draw=black, fill=gray!75, circle, minimum size=3pt, inner sep=0pt, line width=0.4pt] at (3.739,0.642) {};
  \node[draw=black, fill=gray!75, circle, minimum size=3pt, inner sep=0pt, line width=0.4pt] at (1.801,3.803) {};
  \node[draw=black, fill=gray!75, circle, minimum size=3pt, inner sep=0pt, line width=0.4pt] at (1.179,2.645) {};
  \node[draw=black, fill=gray!75, circle, minimum size=3pt, inner sep=0pt, line width=0.4pt] at (2.874,1.837) {};
  \node[draw=black, fill=gray!75, circle, minimum size=3pt, inner sep=0pt, line width=0.4pt] at (3.099,2.625) {};
  \node[draw=black, fill=gray!75, circle, minimum size=3pt, inner sep=0pt, line width=0.4pt] at (3.485,1.888) {};
  \node[draw=black, fill=gray!75, circle, minimum size=3pt, inner sep=0pt, line width=0.4pt] at (2.582,2.336) {};
  \node[draw=black, fill=gray!75, circle, minimum size=3pt, inner sep=0pt, line width=0.4pt] at (2.199,2.598) {};
  \node[draw=black, fill=gray!75, circle, minimum size=3pt, inner sep=0pt, line width=0.4pt] at (2.328,3.024) {};
  \node[draw=black, fill=gray!75, circle, minimum size=3pt, inner sep=0pt, line width=0.4pt] at (3.367,0.761) {};
  \node[draw=black, fill=gray!75, circle, minimum size=3pt, inner sep=0pt, line width=0.4pt] at (3.783,0.413) {};
  \node[draw=black, fill=gray!75, circle, minimum size=3pt, inner sep=0pt, line width=0.4pt] at (0.551,4.011) {};
  \node[draw=black, fill=gray!75, circle, minimum size=3pt, inner sep=0pt, line width=0.4pt] at (1.506,2.919) {};
  \node[draw=black, fill=gray!75, circle, minimum size=3pt, inner sep=0pt, line width=0.4pt] at (2.624,0.553) {};
  \node[draw=black, fill=gray!75, circle, minimum size=3pt, inner sep=0pt, line width=0.4pt] at (1.362,3.562) {};
  \node[draw=black, fill=white, rectangle, minimum size=6pt, inner sep=0pt, line width=0.6pt] at (3.589,2.759) {};
\end{tikzpicture}
\caption{\texttt{s17\_t2366}}
\end{subfigure}
\hspace{1cm}
\begin{subfigure}[b]{0.45\textwidth}
\centering
\begin{tikzpicture}[x=1cm, y=1cm]
  \draw[black, line width=0.8pt] (2.764,1.687) -- (2.591,1.394);
  \draw[black, line width=0.8pt] (1.948,2.211) -- (2.359,2.263);
  \draw[black, line width=0.8pt] (3.416,1.952) -- (3.488,1.329);
  \draw[black, line width=0.8pt] (1.948,2.211) -- (1.214,2.680);
  \draw[black, line width=0.8pt] (1.214,2.680) -- (1.780,3.282);
  \draw[black, line width=0.8pt] (3.416,1.952) -- (4.151,2.544);
  \draw[black, line width=0.8pt] (4.151,2.544) -- (4.746,2.576);
  \draw[black, line width=0.8pt] (5.497,2.293) -- (6.000,2.038);
  \draw[black, line width=0.8pt] (2.591,1.394) -- (2.330,0.180);
  \draw[black, line width=0.8pt] (2.617,0.000) -- (2.330,0.180);
  \draw[black, line width=0.8pt] (2.764,1.687) -- (3.488,1.329);
  \draw[black, line width=0.8pt] (0.889,1.490) -- (0.975,1.283);
  \draw[black, line width=0.8pt] (2.591,1.394) -- (0.975,1.283);
  \draw[black, line width=0.8pt] (1.214,2.680) -- (1.177,2.227);
  \draw[black, line width=0.8pt] (0.889,1.490) -- (1.177,2.227);
  \draw[black, line width=0.8pt] (4.746,2.576) -- (4.221,2.164);
  \draw[black, line width=0.8pt] (2.764,1.687) -- (2.349,2.007);
  \draw[black, line width=0.8pt] (2.359,2.263) -- (2.349,2.007);
  \draw[black, line width=0.8pt] (2.764,1.687) -- (3.118,1.970);
  \draw[black, line width=0.8pt] (3.416,1.952) -- (3.118,1.970);
  \draw[black, line width=0.8pt] (2.617,0.000) -- (2.935,0.063);
  \draw[black, line width=0.8pt] (5.497,2.293) -- (4.892,2.676);
  \draw[black, line width=0.8pt] (4.746,2.576) -- (4.892,2.676);
  \draw[black, line width=0.8pt] (3.488,1.329) -- (3.831,1.901);
  \draw[black, line width=0.8pt] (4.221,2.164) -- (3.831,1.901);
  \draw[black, line width=0.8pt] (3.488,1.329) -- (3.680,0.480);
  \draw[black, line width=0.8pt] (3.839,0.344) -- (3.680,0.480);
  \draw[black, line width=0.8pt] (2.511,3.675) -- (2.194,3.311);
  \draw[black, line width=0.8pt] (1.780,3.282) -- (2.194,3.311);
  \draw[black, line width=0.8pt] (2.359,2.263) -- (2.269,2.889);
  \draw[black, line width=0.8pt] (2.511,3.675) -- (2.269,2.889);
  \draw[black, line width=0.8pt] (1.948,2.211) -- (2.385,1.510);
  \draw[black, line width=0.8pt] (2.591,1.394) -- (2.385,1.510);
  \draw[black, line width=0.8pt] (2.511,3.675) -- (3.477,3.381);
  \draw[black, line width=0.8pt] (3.783,3.630) -- (3.477,3.381);
  \draw[black, line width=0.8pt] (2.359,2.263) -- (2.115,2.734);
  \draw[black, line width=0.8pt] (1.780,3.282) -- (2.115,2.734);
  \draw[black, line width=0.8pt] (4.151,2.544) -- (3.769,3.342);
  \draw[black, line width=0.8pt] (3.783,3.630) -- (3.769,3.342);
  \draw[black, line width=0.8pt] (1.214,2.680) -- (1.008,2.525);
  \draw[black, line width=0.8pt] (0.000,2.795) -- (1.008,2.525);
  \draw[black, line width=0.8pt] (2.359,2.263) -- (2.774,2.307);
  \draw[black, line width=0.8pt] (4.151,2.544) -- (2.828,2.082);
  \draw[black, line width=0.8pt] (2.774,2.307) -- (2.828,2.082);
  \draw[black, line width=0.8pt] (3.488,1.329) -- (3.103,0.407);
  \draw[black, line width=0.8pt] (2.935,0.063) -- (3.103,0.407);
  \node[draw=black, fill=white, circle, minimum size=4.5pt, inner sep=0pt, line width=0.4pt] at (2.764,1.687) {};
  \node[draw=black, fill=white, circle, minimum size=4.5pt, inner sep=0pt, line width=0.4pt] at (1.948,2.211) {};
  \node[draw=black, fill=white, circle, minimum size=4.5pt, inner sep=0pt, line width=0.4pt] at (3.416,1.952) {};
  \node[draw=black, fill=white, circle, minimum size=4.5pt, inner sep=0pt, line width=0.4pt] at (1.214,2.680) {};
  \node[draw=black, fill=white, circle, minimum size=4.5pt, inner sep=0pt, line width=0.4pt] at (0.889,1.490) {};
  \node[draw=black, fill=white, circle, minimum size=4.5pt, inner sep=0pt, line width=0.4pt] at (4.151,2.544) {};
  \node[draw=black, fill=white, circle, minimum size=4.5pt, inner sep=0pt, line width=0.4pt] at (2.511,3.675) {};
  \node[draw=black, fill=white, circle, minimum size=4.5pt, inner sep=0pt, line width=0.4pt] at (2.591,1.394) {};
  \node[draw=black, fill=white, circle, minimum size=4.5pt, inner sep=0pt, line width=0.4pt] at (1.780,3.282) {};
  \node[draw=black, fill=white, circle, minimum size=4.5pt, inner sep=0pt, line width=0.4pt] at (5.497,2.293) {};
  \node[draw=black, fill=white, circle, minimum size=4.5pt, inner sep=0pt, line width=0.4pt] at (3.488,1.329) {};
  \node[draw=black, fill=white, circle, minimum size=4.5pt, inner sep=0pt, line width=0.4pt] at (4.746,2.576) {};
  \node[draw=black, fill=white, circle, minimum size=4.5pt, inner sep=0pt, line width=0.4pt] at (6.000,2.038) {};
  \node[draw=black, fill=white, circle, minimum size=4.5pt, inner sep=0pt, line width=0.4pt] at (0.000,2.795) {};
  \node[draw=black, fill=white, circle, minimum size=4.5pt, inner sep=0pt, line width=0.4pt] at (3.783,3.630) {};
  \node[draw=black, fill=white, circle, minimum size=4.5pt, inner sep=0pt, line width=0.4pt] at (2.617,0.000) {};
  \node[draw=black, fill=white, circle, minimum size=4.5pt, inner sep=0pt, line width=0.4pt] at (3.839,0.344) {};
  \node[draw=black, fill=white, circle, minimum size=4.5pt, inner sep=0pt, line width=0.4pt] at (2.330,0.180) {};
  \node[draw=black, fill=gray!75, circle, minimum size=3pt, inner sep=0pt, line width=0.4pt] at (0.975,1.283) {};
  \node[draw=black, fill=gray!75, circle, minimum size=3pt, inner sep=0pt, line width=0.4pt] at (1.177,2.227) {};
  \node[draw=black, fill=gray!75, circle, minimum size=3pt, inner sep=0pt, line width=0.4pt] at (4.221,2.164) {};
  \node[draw=black, fill=gray!75, circle, minimum size=3pt, inner sep=0pt, line width=0.4pt] at (2.349,2.007) {};
  \node[draw=black, fill=gray!75, circle, minimum size=3pt, inner sep=0pt, line width=0.4pt] at (3.118,1.970) {};
  \node[draw=black, fill=gray!75, circle, minimum size=3pt, inner sep=0pt, line width=0.4pt] at (2.935,0.063) {};
  \node[draw=black, fill=gray!75, circle, minimum size=3pt, inner sep=0pt, line width=0.4pt] at (4.892,2.676) {};
  \node[draw=black, fill=gray!75, circle, minimum size=3pt, inner sep=0pt, line width=0.4pt] at (3.831,1.901) {};
  \node[draw=black, fill=gray!75, circle, minimum size=3pt, inner sep=0pt, line width=0.4pt] at (3.680,0.480) {};
  \node[draw=black, fill=gray!75, circle, minimum size=3pt, inner sep=0pt, line width=0.4pt] at (2.194,3.311) {};
  \node[draw=black, fill=gray!75, circle, minimum size=3pt, inner sep=0pt, line width=0.4pt] at (2.269,2.889) {};
  \node[draw=black, fill=gray!75, circle, minimum size=3pt, inner sep=0pt, line width=0.4pt] at (2.385,1.510) {};
  \node[draw=black, fill=gray!75, circle, minimum size=3pt, inner sep=0pt, line width=0.4pt] at (3.477,3.381) {};
  \node[draw=black, fill=gray!75, circle, minimum size=3pt, inner sep=0pt, line width=0.4pt] at (2.115,2.734) {};
  \node[draw=black, fill=gray!75, circle, minimum size=3pt, inner sep=0pt, line width=0.4pt] at (3.769,3.342) {};
  \node[draw=black, fill=gray!75, circle, minimum size=3pt, inner sep=0pt, line width=0.4pt] at (1.008,2.525) {};
  \node[draw=black, fill=gray!75, circle, minimum size=3pt, inner sep=0pt, line width=0.4pt] at (2.774,2.307) {};
  \node[draw=black, fill=gray!75, circle, minimum size=3pt, inner sep=0pt, line width=0.4pt] at (2.828,2.082) {};
  \node[draw=black, fill=gray!75, circle, minimum size=3pt, inner sep=0pt, line width=0.4pt] at (3.103,0.407) {};
  \node[draw=black, fill=white, rectangle, minimum size=6pt, inner sep=0pt, line width=0.6pt] at (2.359,2.263) {};
\end{tikzpicture}
\caption{\texttt{s19\_t2744}}
\end{subfigure}
\caption{Synthetic infrastructure networks. Small stations are shown as grey circles, large stations as white circles, and the crew base of the single-depot variant as a white square. Instances are named \texttt{sS\_tT}, where \texttt{S} denotes the number of large stations and \texttt{T} the number of tasks.}
\label{fig:instances}
\end{figure}

%% file: figs/validation.tex
\begin{figure}[htbp]
\centering
\begin{subfigure}[b]{0.48\textwidth}
\centering
\begin{tikzpicture}
\begin{axis}[
    width=\textwidth,
    xlabel={Tasks},
    ylabel={Average task duration (min)},
    tick label style={font=\tiny},
    label style={font=\small},
    mark size=1.5pt
]
\addplot[only marks, mark=*, mark options={fill=gray!75, draw=black}] coordinates {
    (1014,30.528599605522682)
    (1016,29.75688976377953)
    (1088,31.30330882352941)
    (928,30.760775862068964)
    (1078,32.13358070500928)
    (1151,31.328410078192874)
    (1204,29.3812292358804)
    (961,29.559833506763788)
    (1049,31.12011439466158)
    (1459,26.473612063056887)
    (1482,31.181511470985154)
    (1498,28.713618157543394)
    (1286,31.845256609642302)
    (1351,25.820133234641006)
    (1372,30.963556851311957)
    (1532,30.595953002610965)
    (1591,28.133249528598366)
    (1598,27.33979974968711)
    (1500,31.832666666666668)
    (1511,28.2819324950364)
    (1629,29.204419889502763)
    (1755,30.022792022792025)
    (1713,30.761821366024517)
    (1745,29.83381088825215)
    (1749,26.252715837621498)
    (2029,28.21981271562346)
    (2052,29.920077972709553)
    (2136,29.854400749063668)
    (2219,26.950878774222623)
    (2252,28.319715808170518)
    (1955,29.46291560102302)
    (2090,29.551674641148328)
    (2166,30.195752539242843)
    (2183,34.36875858909757)
    (2214,28.77416440831075)
    (2366,29.216398985629755)
    (2395,33.3678496868476)
    (2550,30.82549019607843)
    (2641,30.30443014009845)
    (2648,26.128776435045317)
    (2315,31.97365010799136)
    (2484,29.315217391304348)
    (2702,27.424500370096226)
    (2711,31.6672814459609)
    (2744,28.274052478134113)
    (2921,31.83224922971585)
    (3016,27.532493368700266)
    (118,59.0)
    (120,41.0)
    (122,22.0)
    (131,42.0)
    (90,42.0)
    (168,46.0)
    (172,38.5)
    (228,46.425438596491226)
    (232,38.01724137931034)
    (262,39.767175572519086)
    (300,35.28)
    (337,43.55786350148368)
    (273,20.23076923076923)
    (290,32.86896551724138)
    (297,29.946127946127945)
    (311,40.03215434083601)
    (313,35.95527156549521)
    (320,42.15625)
    (333,19.45945945945946)
    (349,17.120343839541547)
    (385,20.10909090909091)
    (409,19.792176039119802)
    (413,19.612590799031477)
    (501,19.756487025948104)
    (591,19.835871404399324)
    (450,27.96888888888889)
    (482,29.42116182572614)
    (514,35.90466926070039)
    (518,27.16023166023166)
    (521,25.809980806142036)
    (580,45.86206896551724)
    (570,34.329824561403505)
    (599,32.70283806343907)
    (627,31.56937799043062)
    (681,24.12481644640235)
    (695,28.359712230215827)
    (738,29.90921409214092)
    (628,34.30891719745223)
    (681,30.93979441997063)
    (694,33.39048991354467)
    (734,28.07356948228883)
    (872,33.934633027522935)
    (773,29.494178525226392)
    (775,32.05290322580645)
    (784,27.184948979591834)
    (946,32.88689217758985)
    (972,27.685185185185183)
    (991,26.507568113017154)
};
\addplot[only marks, mark=*, mark options={fill=white, draw=black}] coordinates {
    (240,29.4875)
    (1125,29.977777777777778)
    (738,28.806233062330623)
    (990,27.311111111111114)
    (1880,25.775)
    (1627,29.255685310387214)
    (184,31.146739130434785)
    (700,26.694285714285716)
    (1130,28.367256637168143)
    (928,27.461206896551726)
    (1064,28.920112781954888)
    (1235,28.38947368421053)
    (287,29.243902439024392)
    (560,24.830357142857146)
    (1942,27.910401647785786)
    (1629,29.990178023327193)
    (1463,27.369104579630896)
    (1661,28.54244431065623)
    (327,28.501529051987767)
    (494,24.87854251012146)
    (1460,25.776712328767125)
    (2196,28.023224043715846)
    (1803,28.509151414309486)
    (2368,26.814611486486488)
    (256,28.02734375)
    (536,23.867537313432834)
    (1144,29.732517482517483)
    (1514,27.41017173051519)
    (2834,28.638673253352156)
    (2923,27.503250085528567)
    (142,39.2112676056338)
    (682,29.146627565982403)
    (672,30.651785714285715)
    (881,31.116912599318955)
    (1107,29.824751580849142)
    (1075,28.57860465116279)
    (357,28.126050420168067)
    (739,27.848443843031124)
    (802,30.452618453865338)
    (1071,31.03174603174603)
    (931,30.408163265306122)
    (1574,31.32909783989835)
};
\end{axis}
\end{tikzpicture}
\caption{Average task duration (min)}
\end{subfigure}
\hfill
\begin{subfigure}[b]{0.48\textwidth}
\centering
\begin{tikzpicture}
\begin{axis}[
    width=\textwidth,
    xlabel={Tasks},
    ylabel={Connectivity (\%)},
    tick label style={font=\tiny},
    label style={font=\small},
    mark size=1.5pt
]
\addplot[only marks, mark=*, mark options={fill=gray!75, draw=black}] coordinates {
    (1014,3.119213537620402)
    (1016,2.945774019626857)
    (1088,2.783903349748363)
    (928,2.9581891901945467)
    (1078,2.7185044693998135)
    (1151,2.608393457484985)
    (1204,2.6491081266932337)
    (961,2.9580731876517516)
    (1049,2.2151230906934267)
    (1459,2.6543538944219267)
    (1482,2.521639370852207)
    (1498,2.178188152004944)
    (1286,2.2229214951800595)
    (1351,2.2105984593031223)
    (1372,2.0769670794232042)
    (1532,2.1100903349915496)
    (1591,2.501413216639193)
    (1598,2.391530427436299)
    (1500,2.1909717589504116)
    (1511,2.294257125450888)
    (1629,2.1729162613140516)
    (1755,2.0458894119099362)
    (1713,1.8650670245674912)
    (1745,2.298769748429326)
    (1749,2.3443275202698373)
    (2029,1.9955711220828558)
    (2052,1.7162264782167782)
    (2136,2.0405406590707753)
    (2219,1.8433310807433627)
    (2252,1.9013850564146346)
    (1955,1.8412228048177128)
    (2090,1.863348915829327)
    (2166,1.7693132795523512)
    (2183,1.8535025883283585)
    (2214,1.644507633508328)
    (2366,1.840949747926492)
    (2395,1.7930002459175076)
    (2550,1.8920453234255645)
    (2641,1.4723819030899683)
    (2648,1.6503035414885687)
    (2315,1.5434270876307423)
    (2484,1.5978217093627973)
    (2702,1.8017561278261114)
    (2711,1.9927288169967645)
    (2744,1.6231084903103474)
    (2921,1.6020503393002021)
    (3016,1.5422005797713465)
    (118,10.343328987396784)
    (120,10.392156862745098)
    (122,10.432190760059612)
    (131,11.180270111567822)
    (90,11.835205992509364)
    (168,8.276304533789563)
    (172,8.32313341493268)
    (228,8.667594095370585)
    (232,8.646812957157785)
    (262,8.212687549355094)
    (300,7.594202898550724)
    (337,7.391903348876642)
    (273,6.135531135531135)
    (290,7.056437179334209)
    (297,6.094731094731095)
    (311,6.298101856653875)
    (313,5.613582370770869)
    (320,6.4498432601880875)
    (333,4.839176525923514)
    (349,4.953397226887988)
    (385,4.847132034632034)
    (409,4.8875785032839545)
    (413,5.600742847739721)
    (501,4.950099800399202)
    (591,5.3411339585305)
    (450,4.52264291017075)
    (482,4.880910275101146)
    (514,4.381034731229284)
    (518,4.167195656557358)
    (521,4.363649785914661)
    (580,4.507772020725389)
    (570,4.160577189899177)
    (599,3.816282432817237)
    (627,4.200233374607008)
    (681,3.772998186058565)
    (695,4.21081002633052)
    (738,3.656881887679121)
    (628,3.350806082954926)
    (681,3.4300768765656042)
    (694,3.9925812260106204)
    (734,3.576433677433265)
    (872,3.155447181874677)
    (773,3.173491343195544)
    (775,2.8943902642327246)
    (784,3.165154950869237)
    (946,3.341499155452644)
    (972,3.163129945370476)
    (991,3.025002803004821)
};
\addplot[only marks, mark=*, mark options={fill=white, draw=black}] coordinates {
    (240,3.985355648535565)
    (1125,2.7297746144721233)
    (738,2.437553547855696)
    (990,2.7496399791647517)
    (1880,2.356221620825926)
    (1627,2.189905734336999)
    (184,4.621050130672369)
    (700,2.8277130594727162)
    (1130,2.236923583404532)
    (928,2.4643826953836996)
    (1064,2.8464269799616635)
    (1235,2.077178984114069)
    (287,2.6680636436733995)
    (560,2.180552006133402)
    (1942,2.3554804954181305)
    (1629,2.66838913247753)
    (1463,2.3029530049473887)
    (1661,2.123557444709603)
    (327,2.34517176037973)
    (494,2.5211257195884076)
    (1460,2.171218793130968)
    (2196,2.0689927015779364)
    (1803,2.269925017066758)
    (2368,2.172324415670423)
    (256,3.4865196078431375)
    (536,2.7346910308271726)
    (1144,1.920170817808613)
    (1514,2.0796426566411226)
    (2834,2.2385131780624614)
    (2923,2.020745565569208)
    (142,5.85356108280891)
    (682,2.697861089221044)
    (672,2.0150273224043715)
    (881,2.8242699411825405)
    (1107,2.259009329092688)
    (1075,2.264778485124074)
    (357,3.5973940137852893)
    (739,3.0617072070585385)
    (802,3.403164996372988)
    (1071,2.1517142682618218)
    (931,2.1754847949366503)
    (1574,2.7107696508181665)
};
\end{axis}
\end{tikzpicture}
\caption{Connectivity (\%)}
\end{subfigure}
\\[1em]
\begin{subfigure}[b]{0.48\textwidth}
\centering
\begin{tikzpicture}
\begin{axis}[
    width=\textwidth,
    xlabel={Tasks},
    ylabel={LB per task},
    tick label style={font=\tiny},
    label style={font=\small},
    mark size=1.5pt
]
\addplot[only marks, mark=*, mark options={fill=gray!75, draw=black}] coordinates {
    (1014,5717.373353648915)
    (1016,5592.2377595472435)
    (1088,5787.129829503676)
    (928,6293.798304418104)
    (1078,6194.494379128015)
    (1151,5873.714674543875)
    (1204,5407.302555481728)
    (961,6236.511341311135)
    (1049,6478.203773307912)
    (1459,5178.009472583962)
    (1482,5717.327697165992)
    (1498,5277.815143791722)
    (1286,6204.124427993779)
    (1351,4972.2001585492235)
    (1372,6087.84779074344)
    (1532,5925.171308028721)
    (1591,5472.287517976115)
    (1598,5330.149208635795)
    (1500,5870.3889188)
    (1511,6016.363800198544)
    (1629,6237.1628834868025)
    (1755,5866.332862849003)
    (1713,5784.706265265616)
    (1745,5951.111721948424)
    (1749,5307.811217724415)
    (2029,5445.835627698374)
    (2052,5827.51288128655)
    (2136,5920.2458895131085)
    (2219,5949.548487336638)
    (2252,5500.444887699822)
    (1955,5958.516189411765)
    (2090,5957.41175861244)
    (2166,6534.124963665743)
    (2183,6599.225148007329)
    (2214,5621.633354923216)
    (2366,5678.629403254437)
    (2395,6410.439179707724)
    (2550,5828.402605058824)
    (2641,5808.051486141613)
    (2648,5274.902360536254)
    (2315,6322.793923585313)
    (2484,6035.886446457327)
    (2702,5303.509701739452)
    (2711,5839.661355256363)
    (2744,5497.449411953353)
    (2921,6145.567155597398)
    (3016,5301.733984880637)
    (118,10276.635094915255)
    (120,8314.25)
    (122,5893.44262295082)
    (131,8697.709923664122)
    (90,8285.333333333334)
    (168,8583.214285714286)
    (172,7044.709302325581)
    (228,8112.566564035087)
    (232,6759.59939137931)
    (262,7063.465791603054)
    (300,6640.305031333333)
    (337,7544.57230860534)
    (273,4255.7392289377285)
    (290,5765.3539210344825)
    (297,5763.5710767676765)
    (311,7299.254383279743)
    (313,6667.094240894568)
    (320,7429.1019559375)
    (333,4983.423423423424)
    (349,3797.9256521489974)
    (385,4267.329649090909)
    (409,4532.460880195599)
    (413,4076.664485230024)
    (501,3854.9870095808383)
    (591,4127.649309137056)
    (450,5440.175238000001)
    (482,5973.101512033195)
    (514,6838.032937743191)
    (518,5265.188621428571)
    (521,5510.995497696737)
    (580,7795.836982931034)
    (570,6240.6812510526315)
    (599,5803.372760267112)
    (627,6037.611925518341)
    (681,4776.857307488987)
    (695,5546.8206297841725)
    (738,5327.78031707317)
    (628,6560.996539171974)
    (681,6312.254639500734)
    (694,6527.053080691643)
    (734,5189.544553814714)
    (872,5954.470282683486)
    (773,5657.441884864165)
    (775,5903.877655741935)
    (784,5431.442149362245)
    (946,5872.48565729387)
    (972,5185.491573456789)
    (991,4839.529297578204)
};
\addplot[only marks, mark=*, mark options={fill=white, draw=black}] coordinates {
    (240,5897.416666666667)
    (1125,6267.114871822222)
    (738,7636.134042140921)
    (990,7665.437174444444)
    (1880,5869.312605638298)
    (1627,6647.918780147511)
    (184,6123.913043478261)
    (700,7253.685434714285)
    (1130,7180.376462035398)
    (928,6926.907145797413)
    (1064,7701.54112518797)
    (1235,9765.82447174089)
    (287,5994.564459930313)
    (560,7865.214285714285)
    (1942,6368.45028753862)
    (1629,6219.381543646409)
    (1463,6708.787472248804)
    (1661,7427.913260024083)
    (327,6291.926605504587)
    (494,6778.686740890688)
    (1460,7125.560224452055)
    (2196,6762.169229098361)
    (1803,7838.595797282307)
    (2368,6356.984822635135)
    (256,5345.859375)
    (536,5239.701492537313)
    (1144,6679.108328758742)
    (1514,6190.407963738441)
    (2834,6252.557659456598)
    (2923,6705.793674751967)
    (142,7291.69014084507)
    (682,6297.731875219941)
    (672,8470.073660714286)
    (881,7142.1317886492625)
    (1107,9214.73755691057)
    (1075,8738.240855813954)
    (357,5390.8403361344535)
    (739,5309.340163193505)
    (802,7185.64253478803)
    (1071,7512.714876937442)
    (931,8859.932189688507)
    (1574,6428.706378780177)
};
\end{axis}
\end{tikzpicture}
\caption{LB per task}
\end{subfigure}
\hfill
\begin{subfigure}[b]{0.48\textwidth}
\centering
\begin{tikzpicture}
\begin{axis}[
    width=\textwidth,
    xlabel={Tasks},
    ylabel={Root node solve time (s)},
    tick label style={font=\tiny},
    label style={font=\small},
    mark size=1.5pt
]
\addplot[only marks, mark=*, mark options={fill=gray!75, draw=black}] coordinates {
    (1014,16.998)
    (1016,20.153)
    (1088,16.15)
    (928,13.056)
    (1078,13.957)
    (1151,19.266)
    (1204,28.907)
    (961,13.687)
    (1049,13.313)
    (1459,44.914)
    (1482,36.729)
    (1498,47.937)
    (1286,19.324)
    (1351,35.039)
    (1372,23.521)
    (1532,36.761)
    (1591,58.565)
    (1598,52.75)
    (1500,32.628)
    (1511,33.888)
    (1629,40.442)
    (1755,69.326)
    (1713,41.783)
    (1745,49.093)
    (1749,77.551)
    (2029,86.22)
    (2052,72.309)
    (2136,137.399)
    (2219,98.548)
    (2252,146.224)
    (1955,55.128)
    (2090,58.809)
    (2166,63.052)
    (2183,78.942)
    (2214,115.672)
    (2366,150.857)
    (2395,106.722)
    (2550,161.131)
    (2641,200.395)
    (2648,535.562)
    (2315,79.083)
    (2484,102.552)
    (2702,252.458)
    (2711,231.827)
    (2744,196.499)
    (2921,194.048)
    (3016,332.628)
    (118,0.156)
    (120,0.114)
    (122,0.158)
    (131,0.109)
    (90,0.093)
    (168,0.17)
    (172,0.239)
    (228,0.389)
    (232,0.603)
    (262,0.431)
    (300,1.157)
    (337,1.104)
    (273,1.731)
    (290,1.051)
    (297,1.211)
    (311,0.964)
    (313,0.774)
    (320,0.786)
    (333,3.061)
    (349,2.809)
    (385,2.136)
    (409,7.755)
    (413,3.844)
    (501,3.458)
    (591,5.76)
    (450,4.521)
    (482,5.139)
    (514,1.892)
    (518,6.874)
    (521,5.792)
    (580,2.755)
    (570,5.761)
    (599,6.766)
    (627,6.251)
    (681,10.846)
    (695,9.065)
    (738,9.506)
    (628,5.668)
    (681,6.449)
    (694,8.299)
    (734,11.096)
    (872,13.695)
    (773,8.412)
    (775,8.755)
    (784,10.053)
    (946,17.394)
    (972,17.097)
    (991,19.968)
};
\addplot[only marks, mark=*, mark options={fill=white, draw=black}] coordinates {
    (240,0.638)
    (1125,12.714)
    (738,2.092)
    (990,6.558)
    (1880,39.154)
    (1627,20.154)
    (184,0.061)
    (700,1.66)
    (1130,7.54)
    (928,5.871)
    (1064,6.921)
    (1235,5.559)
    (287,0.096)
    (560,0.539)
    (1942,50.83)
    (1629,38.611)
    (1463,16.992)
    (1661,14.435)
    (327,0.124)
    (494,0.887)
    (1460,11.882)
    (2196,55.96)
    (1803,22.841)
    (2368,67.882)
    (256,0.09)
    (536,1.703)
    (1144,7.242)
    (1514,13.779)
    (2834,133.308)
    (2923,135.341)
    (142,0.085)
    (682,3.539)
    (672,1.817)
    (881,5.967)
    (1107,4.34)
    (1075,6.429)
    (357,1.447)
    (739,5.151)
    (802,5.152)
    (1071,6.412)
    (931,2.936)
    (1574,20.724)
};
\end{axis}
\end{tikzpicture}
\caption{Root node solve time (s)}
\end{subfigure}
\\[1em]
\begin{subfigure}[b]{0.48\textwidth}
\centering
\begin{tikzpicture}
\begin{axis}[
    width=\textwidth,
    xlabel={Tasks},
    ylabel={Average number of duties},
    tick label style={font=\tiny},
    label style={font=\small},
    mark size=1.5pt
]
\addplot[only marks, mark=*, mark options={fill=gray!75, draw=black}] coordinates {
    (1014,100.4532)
    (1016,98.1045)
    (1088,108.0962)
    (928,100.6566)
    (1078,113.3)
    (1151,114.1778)
    (1204,110.8909)
    (961,104.7225)
    (1049,115.5282)
    (1459,128.7469)
    (1482,140.9099)
    (1498,131.9884)
    (1286,135.4661)
    (1351,115.8986)
    (1372,143.0407)
    (1532,151.9723)
    (1591,152.0)
    (1598,145.53)
    (1500,151.1686)
    (1511,152.607)
    (1629,176.0)
    (1755,174.0)
    (1713,166.3471)
    (1745,173.941)
    (1749,160.0)
    (2029,184.8318)
    (2052,199.092)
    (2136,213.413)
    (2219,223.0498)
    (2252,209.418)
    (1955,197.9863)
    (2090,208.4793)
    (2166,237.9605)
    (2183,241.7697)
    (2214,216.0117)
    (2366,227.3411)
    (2395,257.133)
    (2550,251.2751)
    (2641,256.4595)
    (2648,236.0)
    (2315,244.257)
    (2484,250.7651)
    (2702,242.7696)
    (2711,265.7155)
    (2744,253.747)
    (2921,303.7554)
    (3016,269.1588)
    (118,21.6373)
    (120,18.5)
    (122,13.3333)
    (131,20.6667)
    (90,14.0)
    (168,26.0)
    (172,21.25)
    (228,33.1069)
    (232,27.3529)
    (262,31.6374)
    (300,35.1447)
    (337,42.8971)
    (273,20.5683)
    (290,29.3525)
    (297,30.2988)
    (311,40.0)
    (313,36.5969)
    (320,41.7445)
    (333,30.0)
    (349,23.5461)
    (385,28.9043)
    (409,34.0)
    (413,29.2162)
    (501,33.3835)
    (591,44.0)
    (450,43.0)
    (482,50.0)
    (514,60.4249)
    (518,46.0)
    (521,49.39)
    (580,78.0579)
    (570,61.4774)
    (599,59.7767)
    (627,65.6074)
    (681,55.9285)
    (695,67.0)
    (738,66.3244)
    (628,69.8623)
    (681,73.4891)
    (694,78.245)
    (734,66.0342)
    (872,88.1729)
    (773,75.9889)
    (775,79.2588)
    (784,73.7701)
    (946,96.3717)
    (972,87.3148)
    (991,80.8088)
};
\addplot[only marks, mark=*, mark options={fill=white, draw=black}] coordinates {
    (240,24.3333)
    (1125,122.0)
    (738,96.4128)
    (990,126.7302)
    (1880,186.918)
    (1627,180.9158)
    (184,20.0)
    (700,85.1182)
    (1130,137.489)
    (928,107.2836)
    (1064,136.5999)
    (1235,204.0)
    (287,30.0)
    (560,75.0)
    (1942,208.9953)
    (1629,172.4349)
    (1463,163.0248)
    (1661,207.6503)
    (327,36.0)
    (494,58.75)
    (1460,174.8955)
    (2196,259.0263)
    (1803,241.0)
    (2368,258.7115)
    (256,24.0)
    (536,50.0)
    (1144,129.1535)
    (1514,156.3802)
    (2834,304.0)
    (2923,337.0)
    (142,18.0)
    (682,74.0)
    (672,96.4)
    (881,106.0583)
    (1107,171.097)
    (1075,165.0)
    (357,33.0)
    (739,66.3487)
    (802,95.9506)
    (1071,133.3118)
    (931,138.7474)
    (1574,170.7062)
};
\end{axis}
\end{tikzpicture}
\caption{Average number of duties}
\end{subfigure}
\hfill
\begin{subfigure}[b]{0.48\textwidth}
\centering
\begin{tikzpicture}
\begin{axis}[
    width=\textwidth,
    xlabel={Tasks},
    ylabel={Average duty length (h)},
    tick label style={font=\tiny},
    label style={font=\small},
    mark size=1.5pt
]
\addplot[only marks, mark=*, mark options={fill=gray!75, draw=black}] coordinates {
    (1014,8.03127586111111)
    (1016,8.087476027777779)
    (1088,8.18003188888889)
    (928,8.118185555555556)
    (1078,8.37164188888889)
    (1151,8.44765625)
    (1204,8.308300194444445)
    (961,7.897278916666666)
    (1049,8.339526166666667)
    (1459,8.299672972222222)
    (1482,8.703103916666667)
    (1498,8.639017722222222)
    (1286,8.360184305555554)
    (1351,8.099906722222222)
    (1372,8.220182194444446)
    (1532,8.591775861111111)
    (1591,7.910835972222222)
    (1598,8.257776611111112)
    (1500,8.18058225)
    (1511,8.547085916666667)
    (1629,8.035887527777778)
    (1755,8.435846388888889)
    (1713,8.547064166666667)
    (1745,8.583990805555556)
    (1749,8.116947611111112)
    (2029,8.606088305555556)
    (2052,8.684161805555554)
    (2136,8.459535527777778)
    (2219,8.441330277777778)
    (2252,8.430459611111111)
    (1955,8.343580638888888)
    (2090,8.589698138888888)
    (2166,8.521080305555556)
    (2183,8.551691694444445)
    (2214,8.005173138888889)
    (2366,8.416398333333333)
    (2395,8.585670638888889)
    (2550,8.430010666666668)
    (2641,8.614130472222222)
    (2648,8.440609055555555)
    (2315,8.646024805555555)
    (2484,8.608216333333333)
    (2702,8.396526416666667)
    (2711,8.549989916666666)
    (2744,8.513607361111111)
    (2921,8.4159855)
    (3016,8.502062444444444)
    (118,7.5678371944444445)
    (120,6.980630638888889)
    (122,6.979166666666667)
    (131,7.314516138888889)
    (90,6.7952380833333335)
    (168,7.405769222222222)
    (172,7.8390849722222224)
    (228,7.519297333333333)
    (232,7.925842305555555)
    (262,8.24859475)
    (300,7.7451756666666665)
    (337,8.463959222222222)
    (273,7.690534166666667)
    (290,7.822556777777778)
    (297,7.693523166666667)
    (311,7.764361888888889)
    (313,7.839215638888889)
    (320,7.819218194444445)
    (333,7.365555555555556)
    (349,7.6369221666666665)
    (385,7.788927749999999)
    (409,7.145232833333333)
    (413,8.007685777777777)
    (501,8.070364694444445)
    (591,7.400509722222222)
    (450,7.8144628888888885)
    (482,7.9946385)
    (514,8.157555277777778)
    (518,8.46961175)
    (521,8.148282166666666)
    (580,8.090573583333335)
    (570,8.072699083333333)
    (599,8.153741944444445)
    (627,8.027936638888889)
    (681,8.156734333333334)
    (695,7.982754305555555)
    (738,8.467476)
    (628,8.382644777777777)
    (681,8.248236222222221)
    (694,8.081167333333333)
    (734,8.023368722222223)
    (872,8.357687472222223)
    (773,7.986263805555556)
    (775,8.03575163888889)
    (784,8.034233138888888)
    (946,8.012563583333332)
    (972,8.034893)
    (991,8.486018194444444)
};
\addplot[only marks, mark=*, mark options={fill=white, draw=black}] coordinates {
    (240,8.157305944444445)
    (1125,8.053060638888889)
    (738,8.236507444444445)
    (990,8.633725833333333)
    (1880,8.398021222222223)
    (1627,8.60712025)
    (184,7.65)
    (700,8.57034625)
    (1130,8.392872027777777)
    (928,8.643767416666666)
    (1064,8.663506916666666)
    (1235,8.422648722222222)
    (287,7.93)
    (560,8.313037027777778)
    (1942,8.437812944444444)
    (1629,8.320766694444444)
    (1463,8.723691888888888)
    (1661,8.504460027777778)
    (327,7.8754629722222225)
    (494,7.832961)
    (1460,8.523071999999999)
    (2196,7.924728472222222)
    (1803,8.289751305555555)
    (2368,8.162729416666666)
    (256,7.839583333333334)
    (536,7.602666666666666)
    (1144,8.433719)
    (1514,8.647958722222223)
    (2834,8.191290583333334)
    (2923,8.156474527777778)
    (142,7.978703694444444)
    (682,8.122571833333334)
    (672,8.401249138888888)
    (881,8.479981694444444)
    (1107,8.560966138888888)
    (1075,7.814156444444444)
    (357,8.199747472222223)
    (739,8.426678194444444)
    (802,8.683600805555555)
    (1071,8.765446833333334)
    (931,8.514011305555556)
    (1574,8.465558083333333)
};
\end{axis}
\end{tikzpicture}
\caption{Average duty length (h)}
\end{subfigure}
\par\medskip
\begin{tikzpicture}[x=1cm, y=1cm]
  \node[draw=black, fill=gray, circle, inner sep=1.5pt] at (0, 0.08) {};
  \node[font=\small, right] at (0.15, 0.08) {Synthetic};
  \node[draw=black, fill=white, circle, inner sep=1.5pt] at (2.5, 0.08) {};
  \node[font=\small, right] at (2.65, 0.08) {NS};
\end{tikzpicture}
\caption{Comparison of synthetic and NS instances.}
\label{fig:validation}
\end{figure}

%% file: tables/branchingTable.tex
\begin{table}[htbp]
\centering
\caption{Branching rule comparison. Every setting branches on duties (D) and connections (C). The additional branching rule is indicated in the rule name.}
\label{tab:branching}
\begin{adjustbox}{max width=\textwidth}
\begin{tabular}{llrrrrrrrr}
\toprule
 &  & \multicolumn{2}{c}{Tiny} & \multicolumn{2}{c}{Small} & \multicolumn{2}{c}{Medium} & \multicolumn{2}{c}{Large} \\
\cmidrule(lr){3-4}\cmidrule(lr){5-6}\cmidrule(lr){7-8}\cmidrule(lr){9-10}
Strong branching & Rule & $\Delta$ LB (\%) & Nodes & $\Delta$ LB (\%) & Nodes & $\Delta$ LB (\%) & Nodes & $\Delta$ LB (\%) & Nodes \\
\midrule
\multirow{6}{*}{None} & None & 0.537 & 4,603 & 0.130 & 3,012 & 0.039 & 1,233 & 0.012 & 450 \\
 & Connection time & 0.537 & 4,165 & 0.129 & 2,969 & 0.039 & 1,267 & 0.013 & 467 \\
 & Depot connection & 0.541 & 6,910 & 0.136 & 3,412 & 0.044 & 1,308 & 0.014 & 545 \\
 & Cumulative start & \textbf{0.557} & 4,742 & 0.135 & 3,493 & 0.039 & 1,284 & 0.012 & 446 \\
 & Station flow & 0.546 & 5,894 & \textbf{0.151} & 3,463 & \textbf{0.045} & 1,297 & \textbf{0.016} & 407 \\
 & Depot flow & 0.547 & 5,910 & 0.132 & 3,862 & 0.038 & 1,506 & 0.014 & 467 \\
\midrule
\multirow{6}{*}{Semi} & None & 0.599 & 284 & \textbf{0.258} & 642 & \textbf{0.089} & 289 & \textbf{0.033} & 100 \\
 & Connection time & \textbf{0.600} & 269 & 0.253 & 651 & \textbf{0.089} & 294 & 0.031 & 100 \\
 & Depot connection & 0.589 & 443 & 0.221 & 676 & 0.074 & 312 & 0.028 & 108 \\
 & Cumulative start & 0.592 & 664 & 0.217 & 638 & 0.067 & 267 & 0.022 & 96 \\
 & Station flow & \textbf{0.600} & 409 & 0.237 & 656 & 0.075 & 309 & 0.028 & 106 \\
 & Depot flow & 0.598 & 449 & 0.208 & 699 & 0.061 & 324 & 0.023 & 113 \\
\midrule
\multirow{6}{*}{Full} & None & \textbf{0.599} & 177 & \textbf{0.232} & 263 & \textbf{0.080} & 121 & \textbf{0.028} & 49 \\
 & Connection time & 0.596 & 154 & 0.231 & 243 & 0.073 & 116 & 0.027 & 50 \\
 & Depot connection & 0.587 & 271 & 0.207 & 292 & 0.066 & 135 & 0.025 & 54 \\
 & Cumulative start & 0.587 & 417 & 0.208 & 286 & 0.062 & 128 & 0.019 & 53 \\
 & Station flow & 0.594 & 271 & 0.217 & 271 & 0.067 & 122 & 0.024 & 51 \\
 & Depot flow & 0.593 & 299 & 0.193 & 304 & 0.055 & 145 & 0.020 & 57 \\
\bottomrule
\end{tabular}
\end{adjustbox}
\end{table}

%% file: figs/connectionRemovalScatter.tex
\begin{figure}[htbp]
\centering
\begin{tikzpicture}
\begin{axis}[
    width=0.7\textwidth,
    xlabel={Tasks},
    ylabel={Lower bound on alternative optimal solutions},
    xtick={150,300,450},
    ytick={0,2000,4000,6000},
    ymax=6500,
    tick label style={font=\tiny},
    label style={font=\small},
    mark size=1.5pt
]
\addplot[only marks, mark=*, mark options={fill=gray!75, draw=black}] coordinates {
    (118,121)
    (120,2)
    (122,29)
    (131,690)
    (90,15)
    (168,7)
    (172,29)
    (228,44)
    (232,243)
    (262,423)
    (300,269)
    (337,1535)
    (273,1026)
    (290,183)
    (297,184)
    (311,565)
    (313,666)
    (320,1585)
    (349,4009)
    (385,2907)
    (409,3753)
    (413,2941)
    (450,949)
    (482,5415)
};
\end{axis}
\end{tikzpicture}
\caption{Lower bound on the number of alternative optimal solutions at the root node versus instance size for all tiny instances.}
\label{fig:connection_removal_scatter}
\end{figure}

%% file: tables/primalTable.tex
\begin{table}[htbp]
\centering
\caption{Primal heuristic comparison.}
\label{tab:primal}
\begin{adjustbox}{max width=\textwidth}
\begin{tabular}{llrrrrrrrr}
\toprule
 &  & \multicolumn{2}{c}{Tiny} & \multicolumn{2}{c}{Small} & \multicolumn{2}{c}{Medium} & \multicolumn{2}{c}{Large} \\
\cmidrule(lr){3-4}\cmidrule(lr){5-6}\cmidrule(lr){7-8}\cmidrule(lr){9-10}
Heuristic & Variant & Gap (\%) & Time (s) & Gap (\%) & Time (s) & Gap (\%) & Time (s) & Gap (\%) & Time (s) \\
\midrule
\multirow{2}{*}{RMH} & Regular & 1.24 & 324 & 2.86 & 1,673 & 4.55 & 1,801 & 6.57 & 1,800 \\
 & With enumeration & \textbf{0.88} & 502 & 1.77 & 1,554 & 4.86 & 1,807 & 8.57 & 1,804 \\
\midrule
\multirow{4}{*}{Diving} & Default & 1.46 & 4 & 1.34 & 31 & 1.06 & 130 & 1.06 & 584 \\
 & Backtracking (1,1) & 1.12 & 9 & 1.07 & 68 & 1.03 & 307 & 0.98 & 1,267 \\
 & Backtracking (1,2) & 1.14 & 14 & \textbf{0.92} & 97 & \textbf{1.00} & 458 & 0.95 & 1,972 \\
 & Backtracking (2,1) & 1.16 & 15 & 0.93 & 90 & 1.01 & 483 & \textbf{0.89} & 2,027 \\
\bottomrule
\end{tabular}
\end{adjustbox}
\end{table}

%% file: tables/tableSolver.tex
\begin{table}[htbp]
\centering
\caption{Branch and price results. Results are averaged over instances of similar size, with a time limit of three hours.}
\label{tab:solver_comparison}
\begin{adjustbox}{max width=\textwidth}
\begin{tabular}{ll rr r r rrrr}
\toprule
& & & & & & \multicolumn{4}{c}{Time (\%)} \\
\cmidrule(lr){7-10}
Group & Configuration & LB (\%) & UB (\%) & Gap (\%) & Nodes & RMP & Pricing & Branching & Heuristic \\
\midrule
\multirow{2}{*}{Tiny} & B\&P & 100.59 & 100.62 & 0.03 & 606 & 21.2 & 5.6 & 39.4 & 31.3 \\
 & B\&P + E + F & 100.59 & 100.64 & 0.05 & 367 & 10.2 & 3.6 & 24.3 & 59.1 \\
\midrule
\multirow{2}{*}{Small} & B\&P & 100.25 & 100.83 & 0.58 & 1,522 & 28.6 & 2.7 & 27.0 & 40.6 \\
 & B\&P + E + F & 100.25 & 100.92 & 0.66 & 1,578 & 25.7 & 2.7 & 26.1 & 42.7 \\
\midrule
\multirow{2}{*}{Medium} & B\&P & 100.09 & 100.99 & 0.89 & 693 & 42.0 & 4.1 & 29.8 & 23.0 \\
 & B\&P + E + F & 100.09 & 100.93 & 0.83 & 797 & 37.5 & 4.5 & 31.1 & 25.9 \\
\midrule
\multirow{2}{*}{Large} & B\&P & 100.03 & 100.89 & 0.85 & 239 & 47.8 & 4.5 & 22.8 & 24.2 \\
 & B\&P + E + F & 100.03 & 100.90 & 0.86 & 280 & 43.9 & 5.5 & 24.5 & 25.6 \\
\midrule
\multirow{2}{*}{NS} & B\&P & 100.14 & 100.33 & 0.19 & 1,203 & 35.0 & 17.6 & 23.9 & 19.2 \\
 & B\&P + E + F & 100.14 & 100.35 & 0.21 & 1,016 & 28.5 & 20.4 & 21.1 & 22.3 \\
\bottomrule
\end{tabular}
\end{adjustbox}
\end{table}

%% file: tables/tableLocalSearch.tex
\begin{table}[htbp]
\centering
\caption{Postprocessing heuristic results, averaged over instances of similar size.}
\label{tab:local_search}
\begin{adjustbox}{max width=\textwidth}
\begin{tabular}{llrrrrrr}
\toprule
 &  & \multicolumn{2}{c}{Small} & \multicolumn{2}{c}{Medium} & \multicolumn{2}{c}{Large} \\
\cmidrule(lr){3-4}\cmidrule(lr){5-6}\cmidrule(lr){7-8}
Fixed (\%) & Strategy & $\Delta$ UB (\%) & Time (s) & $\Delta$ UB (\%) & Time (s) & $\Delta$ UB (\%) & Time (s) \\
\midrule
\multirow{3}{*}{25} & Random & \textbf{0.17} & 271 & \textbf{0.18} & 1,455 & \textbf{0.16} & 5,512 \\
 & Cost per task & 0.10 & 160 & 0.06 & 969 & 0.10 & 3,837 \\
 & Time window & 0.13 & 519 & 0.16 & 2,374 & 0.12 & 7,326 \\
\midrule
\multirow{3}{*}{50} & Random & 0.11 & 115 & 0.11 & 372 & 0.09 & 1,248 \\
 & Cost per task & 0.02 & 18 & 0.02 & 38 & 0.04 & 347 \\
 & Time window & 0.11 & 160 & 0.13 & 801 & 0.09 & 3,420 \\
\midrule
\multirow{3}{*}{75} & Random & 0.04 & 60 & 0.03 & 137 & 0.03 & 386 \\
 & Cost per task & 0.01 & 3 & 0.00 & 5 & 0.01 & 21 \\
 & Time window & 0.02 & 45 & 0.04 & 124 & 0.04 & 407 \\
\bottomrule
\end{tabular}
\end{adjustbox}
\end{table}

%% file: tables/tableMCF.tex
\begin{table}[htbp]
\centering
\caption{Comparison of branch and price and the multi-commodity flow (MCF) formulation with and without valid inequalities (VI) on the tiny instances.}
\label{tab:mcf}
\begin{adjustbox}{max width=\textwidth}
\begin{tabular}{l rrr rrr rrr}
\toprule
 & \multicolumn{3}{c}{Branch and price} & \multicolumn{3}{c}{MCF} & \multicolumn{3}{c}{MCF + VI} \\
\cmidrule(lr){2-4}\cmidrule(lr){5-7}\cmidrule(lr){8-10}
Instance & LB & Gap (\%) & Time (s) & LB & Gap (\%) & Time (s) & LB & Gap (\%) & Time (s) \\
\midrule
\texttt{s2\_t90} & 748,080 & \textbf{0} & \textbf{3} & 748,080 & \textbf{0} & 6,859 & 748,080 & \textbf{0} & 661 \\
\texttt{s2\_t118} & 1,221,780 & \textbf{0} & \textbf{45} & 1,178,195 & 3.64 & 10,800 & 1,212,939 & 0.72 & 10,800 \\
\texttt{s2\_t120} & 1,010,640 & \textbf{0} & \textbf{15} & 937,746 & 7.30 & 10,800 & 987,337 & 2.31 & 10,800 \\
\texttt{s2\_t122} & 740,280 & \textbf{0} & \textbf{144} & 696,039 & 5.98 & 10,800 & 716,498 & 3.21 & 10,800 \\
\texttt{s2\_t131} & 1,147,500 & \textbf{0} & \textbf{29} & 1,094,649 & 4.68 & 10,800 & 1,144,839 & 0.23 & 10,800 \\
\texttt{s3\_t168} & 1,441,980 & \textbf{0} & \textbf{1} & 1,360,521 & 8.47 & 10,800 & 1,441,980 & \textbf{0} & 5,789 \\
\texttt{s3\_t172} & 1,228,380 & \textbf{0} & \textbf{265} & 1,148,631 & 11.49 & 10,800 & 1,205,889 & 2.59 & 10,800 \\
\texttt{s3\_t228} & 1,855,440 & \textbf{0} & \textbf{19} & 1,764,934 & 13.02 & 10,800 & 1,834,120 & 7.87 & 10,800 \\
\texttt{s3\_t232} & 1,575,840 & \textbf{0} & \textbf{7} & 1,514,221 & 12.78 & 10,800 & 1,570,920 & 0.31 & 10,800 \\
\texttt{s3\_t262} & 1,862,520 & \textbf{0} & \textbf{1,879} & 1,753,868 & -- & 10,800 & 1,830,219 & -- & 10,800 \\
\texttt{s3\_t300} & 2,005,200 & \textbf{0} & \textbf{44} & 1,899,559 & -- & 10,800 & 1,986,632 & -- & 10,800 \\
\texttt{s3\_t337} & 2,552,400 & \textbf{0} & \textbf{5,859} & 2,392,853 & -- & 10,800 & 2,524,396 & -- & 10,800 \\
\texttt{s4\_t273} & 1,169,640 & \textbf{0} & \textbf{847} & 1,092,465 & 10.54 & 10,800 & 1,161,610 & 0.69 & 10,800 \\
\texttt{s4\_t290} & 1,676,820 & \textbf{0} & \textbf{13} & 1,632,837 & -- & 10,800 & 1,676,820 & \textbf{0} & 3,108 \\
\texttt{s4\_t297} & 1,725,720 & \textbf{0} & \textbf{34} & 1,622,711 & -- & 10,800 & 1,698,184 & 8.19 & 10,800 \\
\texttt{s4\_t311} & 2,272,440 & \textbf{0} & \textbf{106} & 2,167,980 & -- & 10,800 & 2,263,073 & -- & 10,800 \\
\texttt{s4\_t313} & 2,095,620 & \textbf{0} & \textbf{416} & 1,968,525 & -- & 10,800 & 2,071,977 & -- & 10,800 \\
\texttt{s4\_t320} & 2,385,300 & \textbf{0} & \textbf{600} & 2,238,235 & -- & 10,800 & 2,360,794 & -- & 10,800 \\
\texttt{s5\_t333} & 1,659,480 & \textbf{0} & \textbf{10} & 1,521,505 & -- & 10,800 & 1,659,480 & \textbf{0} & 3,323 \\
\texttt{s5\_t349} & 1,333,484 & \textbf{0.01} & \textbf{10,800} & 1,258,318 & -- & \textbf{10,800} & 1,308,353 & 13.45 & \textbf{10,800} \\
\texttt{s5\_t385} & 1,646,060 & \textbf{0.43} & \textbf{10,800} & 1,546,342 & -- & \textbf{10,800} & 1,638,550 & -- & \textbf{10,800} \\
\texttt{s5\_t409} & 1,854,900 & \textbf{0} & \textbf{2,180} & 1,711,636 & -- & 10,800 & 1,852,650 & -- & 10,800 \\
\texttt{s5\_t413} & 1,700,778 & \textbf{0.08} & \textbf{10,800} & 1,590,711 & -- & \textbf{10,800} & 1,681,564 & -- & \textbf{10,800} \\
\texttt{s6\_t450} & 2,451,294 & \textbf{0.12} & \textbf{10,800} & 2,268,485 & -- & \textbf{10,800} & 2,434,557 & -- & \textbf{10,800} \\
\texttt{s6\_t482} & 2,863,872 & \textbf{0.15} & \textbf{10,800} & 2,687,079 & -- & \textbf{10,800} & 2,844,474 & -- & \textbf{10,800} \\
\midrule
\textbf{Avg.} & & 0.03 & 2,661 & & 7.79 & 10,642 & & 2.83 & 9,587 \\
\bottomrule
\end{tabular}
\end{adjustbox}
\end{table}

%% file: tables/tableSolverMultipleDepots.tex
\begin{table}[htbp]
\centering
\caption{Branch and price results for multi-depot instances. Results are averaged over instances of similar size, with a time limit of three hours.}
\label{tab:solver_comparison_multi_depot}
\begin{adjustbox}{max width=\textwidth}
\begin{tabular}{ll rr r r rrrr}
\toprule
& & & & & & \multicolumn{4}{c}{Time (\%)} \\
\cmidrule(lr){7-10}
Group & Configuration & LB (\%) & UB (\%) & Gap (\%) & Nodes & RMP & Pricing & Branching & Heuristic \\
\midrule
\multirow{2}{*}{Tiny} & B\&P & 100.36 & 100.37 & 0.01 & 968 & 33.6 & 9.5 & 30.0 & 22.4 \\
 & B\&P + E + F & 100.36 & 100.37 & 0.01 & 417 & 22.0 & 9.1 & 16.3 & 46.5 \\
\midrule
\multirow{2}{*}{Small} & B\&P & 100.11 & 100.42 & 0.31 & 1,750 & 32.4 & 7.5 & 23.7 & 33.3 \\
 & B\&P + E + F & 100.11 & 100.48 & 0.37 & 1,544 & 27.2 & 6.8 & 20.0 & 38.7 \\
\midrule
\multirow{2}{*}{Medium} & B\&P & 100.04 & 100.56 & 0.52 & 623 & 37.2 & 11.9 & 22.0 & 25.9 \\
 & B\&P + E + F & 100.04 & 100.54 & 0.50 & 649 & 36.0 & 12.0 & 22.3 & 26.7 \\
\midrule
\multirow{2}{*}{Large} & B\&P & 100.01 & 100.49 & 0.47 & 109 & 23.6 & 7.8 & 11.4 & 56.2 \\
 & B\&P + E + F & 100.01 & 100.51 & 0.50 & 115 & 22.7 & 7.7 & 11.6 & 57.0 \\
\bottomrule
\end{tabular}
\end{adjustbox}
\end{table}

%% file: tables/tableBestKnownTiny.tex
\begin{table}[htbp]
\centering
\caption{Best known bounds and solutions for tiny instances.}
\label{tab:best_tiny}
\begin{adjustbox}{max width=\textwidth}
\begin{tabular}{l rrr r r rrr r r}
\toprule
& \multicolumn{5}{c}{Single-depot} & \multicolumn{5}{c}{Multi-depot} \\
\cmidrule(lr){2-6} \cmidrule(lr){7-11}
Instance & Root LB & Best LB & Best UB & Gap (\%) & Duties & Root LB & Best LB & Best UB & Gap (\%) & Duties \\
\midrule
\texttt{s2\_t90} & 745,680 & 748,080 & 748,080 & 0 & 14 & 728,800 & 748,080 & 748,080 & 0 & 14 \\
\texttt{s2\_t118} & 1,212,642 & 1,221,780 & 1,221,780 & 0 & 22 & 1,200,240 & 1,201,140 & 1,201,140 & 0 & 21 \\
\texttt{s2\_t120} & 997,710 & 1,010,640 & 1,010,640 & 0 & 19 & 975,560 & 983,880 & 983,880 & 0 & 18 \\
\texttt{s2\_t122} & 719,000 & 740,280 & 740,280 & 0 & 14 & 711,060 & 716,460 & 716,460 & 0 & 13 \\
\texttt{s2\_t131} & 1,139,400 & 1,147,500 & 1,147,500 & 0 & 21 & 1,123,140 & 1,123,920 & 1,123,920 & 0 & 20 \\
\texttt{s3\_t168} & 1,441,980 & 1,441,980 & 1,441,980 & 0 & 26 & 1,411,155 & 1,414,680 & 1,414,680 & 0 & 26 \\
\texttt{s3\_t172} & 1,211,690 & 1,228,380 & 1,228,380 & 0 & 22 & 1,199,043 & 1,209,480 & 1,209,480 & 0 & 21 \\
\texttt{s3\_t228} & 1,849,665 & 1,855,440 & 1,855,440 & 0 & 33 & 1,811,664 & 1,819,200 & 1,819,200 & 0 & 33 \\
\texttt{s3\_t232} & 1,568,227 & 1,575,840 & 1,575,840 & 0 & 27 & 1,558,232 & 1,561,560 & 1,561,560 & 0 & 27 \\
\texttt{s3\_t262} & 1,850,628 & 1,862,520 & 1,862,520 & 0 & 32 & 1,749,558 & 1,761,300 & 1,761,300 & 0 & 31 \\
\texttt{s3\_t300} & 1,992,091 & 2,005,200 & 2,005,200 & 0 & 36 & 1,939,398 & 1,941,660 & 1,941,660 & 0 & 35 \\
\texttt{s3\_t337} & 2,542,520 & 2,552,400 & 2,552,400 & 0 & 43 & 2,375,106 & 2,383,200 & 2,383,200 & 0 & 41 \\
\texttt{s4\_t273} & 1,161,816 & 1,169,640 & 1,169,640 & 0 & 21 & 1,117,997 & 1,118,460 & 1,118,460 & 0 & 20 \\
\texttt{s4\_t290} & 1,671,952 & 1,676,820 & 1,676,820 & 0 & 30 & 1,631,119 & 1,633,680 & 1,633,680 & 0 & 29 \\
\texttt{s4\_t297} & 1,711,780 & 1,725,720 & 1,725,720 & 0 & 31 & 1,660,696 & 1,661,280 & 1,661,280 & 0 & 30 \\
\texttt{s4\_t311} & 2,270,068 & 2,272,440 & 2,272,440 & 0 & 40 & 2,205,150 & 2,208,360 & 2,208,360 & 0 & 40 \\
\texttt{s4\_t313} & 2,086,800 & 2,095,620 & 2,095,620 & 0 & 37 & 1,994,272 & 1,996,980 & 1,996,980 & 0 & 35 \\
\texttt{s4\_t320} & 2,377,312 & 2,385,300 & 2,385,300 & 0 & 42 & 2,313,277 & 2,317,440 & 2,317,440 & 0 & 41 \\
\texttt{s5\_t333} & 1,659,480 & 1,659,480 & 1,659,480 & 0 & 30 & 1,584,592 & 1,584,840 & 1,584,840 & 0 & 30 \\
\texttt{s5\_t349} & 1,325,476 & 1,333,560 & 1,333,560 & 0 & 24 & 1,250,717 & 1,258,080 & 1,258,080 & 0 & 22 \\
\texttt{s5\_t385} & 1,642,921 & 1,651,380 & 1,653,120 & 0.11 & 29 & 1,570,116 & 1,571,880 & 1,571,880 & 0 & 28 \\
\texttt{s5\_t409} & 1,853,776 & 1,854,900 & 1,854,900 & 0 & 34 & 1,813,958 & 1,814,700 & 1,814,700 & 0 & 34 \\
\texttt{s5\_t413} & 1,683,662 & 1,701,840 & 1,701,840 & 0 & 30 & 1,623,060 & 1,624,620 & 1,624,620 & 0 & 29 \\
\texttt{s6\_t450} & 2,444,367 & 2,451,720 & 2,451,720 & 0 & 43 & 2,363,507 & 2,368,593 & 2,370,540 & 0.08 & 42 \\
\texttt{s6\_t482} & 2,861,635 & 2,864,880 & 2,868,240 & 0.12 & 50 & 2,688,349 & 2,691,240 & 2,691,240 & 0 & 49 \\
\bottomrule
\end{tabular}
\end{adjustbox}
\end{table}

%% file: tables/tableBestKnownSmall.tex
\begin{table}[htbp]
\centering
\caption{Best known bounds and solutions for small instances.}
\label{tab:best_small}
\begin{adjustbox}{max width=\textwidth}
\begin{tabular}{l rrr r r rrr r r}
\toprule
& \multicolumn{5}{c}{Single-depot} & \multicolumn{5}{c}{Multi-depot} \\
\cmidrule(lr){2-6} \cmidrule(lr){7-11}
Instance & Root LB & Best LB & Best UB & Gap (\%) & Duties & Root LB & Best LB & Best UB & Gap (\%) & Duties \\
\midrule
\texttt{s5\_t501} & 1,931,348 & 1,942,468 & 1,951,260 & 0.45 & 34 & 1,871,901 & 1,878,394 & 1,879,800 & 0.07 & 33 \\
\texttt{s5\_t591} & 2,439,440 & 2,440,976 & 2,441,040 & 0.00 & 44 & 2,341,780 & 2,343,393 & 2,346,600 & 0.14 & 42 \\
\texttt{s6\_t514} & 3,514,748 & 3,528,422 & 3,532,800 & 0.12 & 61 & 3,196,426 & 3,200,895 & 3,202,320 & 0.04 & 57 \\
\texttt{s6\_t518} & 2,708,987 & 2,712,000 & 2,712,000 & 0 & 46 & 2,558,366 & 2,559,996 & 2,560,920 & 0.04 & 46 \\
\texttt{s6\_t521} & 2,857,952 & 2,869,128 & 2,880,120 & 0.38 & 50 & 2,690,373 & 2,695,216 & 2,696,580 & 0.05 & 48 \\
\texttt{s6\_t580} & 4,521,585 & 4,534,807 & 4,561,380 & 0.58 & 78 & 4,312,151 & 4,320,132 & 4,330,920 & 0.25 & 76 \\
\texttt{s7\_t570} & 3,572,503 & 3,588,419 & 3,589,680 & 0.04 & 63 & 3,375,250 & 3,385,253 & 3,391,680 & 0.19 & 60 \\
\texttt{s7\_t599} & 3,473,001 & 3,482,597 & 3,482,820 & 0.01 & 61 & 3,280,091 & 3,284,272 & 3,297,480 & 0.40 & 58 \\
\texttt{s7\_t627} & 3,785,095 & 3,796,950 & 3,799,800 & 0.08 & 67 & 3,536,373 & 3,540,021 & 3,545,940 & 0.17 & 63 \\
\texttt{s7\_t681} & 3,315,261 & 3,333,440 & 3,358,680 & 0.75 & 58 & 3,084,196 & 3,087,013 & 3,096,240 & 0.30 & 55 \\
\texttt{s7\_t695} & 3,867,536 & 3,873,529 & 3,885,720 & 0.31 & 67 & 3,619,377 & 3,622,008 & 3,625,620 & 0.10 & 65 \\
\texttt{s7\_t738} & 3,958,744 & 3,968,197 & 3,985,920 & 0.44 & 68 & 3,657,380 & 3,661,081 & 3,668,820 & 0.21 & 65 \\
\texttt{s8\_t628} & 4,120,124 & 4,129,320 & 4,173,900 & 1.07 & 72 & 3,736,804 & 3,739,828 & 3,746,820 & 0.19 & 65 \\
\texttt{s8\_t681} & 4,257,209 & 4,266,438 & 4,294,380 & 0.65 & 74 & 3,883,493 & 3,886,140 & 3,886,140 & 0 & 71 \\
\texttt{s8\_t694} & 4,541,909 & 4,553,863 & 4,577,040 & 0.51 & 80 & 4,270,748 & 4,273,749 & 4,286,700 & 0.30 & 77 \\
\texttt{s8\_t734} & 3,793,613 & 3,803,805 & 3,825,540 & 0.57 & 67 & 3,516,953 & 3,519,694 & 3,526,980 & 0.21 & 63 \\
\texttt{s8\_t872} & 5,253,048 & 5,265,660 & 5,323,500 & 1.09 & 92 & 4,904,474 & 4,907,288 & 4,924,980 & 0.36 & 86 \\
\texttt{s9\_t773} & 4,373,174 & 4,384,561 & 4,409,760 & 0.57 & 77 & 4,054,312 & 4,057,377 & 4,064,340 & 0.17 & 73 \\
\texttt{s9\_t775} & 4,566,158 & 4,580,470 & 4,608,780 & 0.61 & 81 & 4,212,572 & 4,216,509 & 4,232,580 & 0.38 & 76 \\
\texttt{s9\_t784} & 4,266,799 & 4,275,539 & 4,284,180 & 0.20 & 75 & 3,982,069 & 3,986,176 & 3,992,460 & 0.16 & 72 \\
\texttt{s9\_t946} & 5,564,229 & 5,572,249 & 5,590,680 & 0.33 & 97 & 5,256,955 & 5,259,692 & 5,272,380 & 0.24 & 95 \\
\texttt{s9\_t972} & 5,076,012 & 5,083,015 & 5,134,260 & 1.00 & 90 & 4,768,967 & 4,771,382 & 4,787,940 & 0.35 & 86 \\
\texttt{s9\_t991} & 4,801,999 & 4,809,359 & 4,831,620 & 0.46 & 82 & 4,480,773 & 4,484,429 & 4,507,920 & 0.52 & 79 \\
\texttt{s10\_t928} & 5,818,100 & 5,831,520 & 5,900,040 & 1.16 & 103 & 5,213,343 & 5,217,440 & 5,235,000 & 0.34 & 96 \\
\texttt{s11\_t961} & 6,036,364 & 6,046,847 & 6,108,540 & 1.01 & 108 & 5,236,197 & 5,239,202 & 5,256,000 & 0.32 & 97 \\
\bottomrule
\end{tabular}
\end{adjustbox}
\end{table}

%% file: tables/tableBestKnownMedium.tex
\begin{table}[htbp]
\centering
\caption{Best known bounds and solutions for medium instances.}
\label{tab:best_medium}
\begin{adjustbox}{max width=\textwidth}
\begin{tabular}{l rrr r r rrr r r}
\toprule
& \multicolumn{5}{c}{Single-depot} & \multicolumn{5}{c}{Multi-depot} \\
\cmidrule(lr){2-6} \cmidrule(lr){7-11}
Instance & Root LB & Best LB & Best UB & Gap (\%) & Duties & Root LB & Best LB & Best UB & Gap (\%) & Duties \\
\midrule
\texttt{s9\_t1020} & 5,656,977 & 5,662,419 & 5,694,780 & 0.57 & 97 & 5,161,683 & 5,165,267 & 5,182,380 & 0.33 & 89 \\
\texttt{s10\_t1014} & 5,784,168 & 5,793,198 & 5,821,860 & 0.49 & 102 & 5,415,092 & 5,418,400 & 5,432,940 & 0.27 & 98 \\
\texttt{s10\_t1016} & 5,722,362 & 5,729,371 & 5,770,740 & 0.72 & 101 & 5,345,895 & 5,348,589 & 5,365,800 & 0.32 & 98 \\
\texttt{s10\_t1088} & 6,352,434 & 6,363,921 & 6,415,740 & 0.81 & 112 & 5,804,036 & 5,806,832 & 5,826,240 & 0.33 & 103 \\
\texttt{s11\_t1078} & 6,694,088 & 6,700,406 & 6,735,120 & 0.52 & 116 & 6,140,054 & 6,142,352 & 6,160,680 & 0.30 & 111 \\
\texttt{s11\_t1151} & 6,781,233 & 6,790,106 & 6,858,120 & 0.99 & 117 & 6,158,369 & 6,161,425 & 6,173,640 & 0.20 & 111 \\
\texttt{s11\_t1204} & 6,532,036 & 6,540,270 & 6,582,900 & 0.65 & 114 & 6,106,999 & 6,109,832 & 6,145,740 & 0.58 & 107 \\
\texttt{s12\_t1049} & 6,878,240 & 6,889,604 & 6,979,860 & 1.29 & 119 & 5,718,836 & 5,722,722 & 5,742,660 & 0.35 & 104 \\
\texttt{s12\_t1459} & 7,570,579 & 7,576,125 & 7,637,520 & 0.80 & 131 & 7,010,843 & 7,012,799 & 7,037,040 & 0.34 & 128 \\
\texttt{s12\_t1482} & 8,496,623 & 8,502,321 & 8,557,440 & 0.64 & 144 & 7,623,706 & 7,628,849 & 7,664,520 & 0.47 & 132 \\
\texttt{s12\_t1498} & 7,951,231 & 7,960,595 & 8,040,960 & 1.00 & 136 & 7,141,899 & 7,145,873 & 7,197,540 & 0.72 & 125 \\
\texttt{s13\_t1286} & 8,070,378 & 8,079,261 & 8,155,980 & 0.94 & 140 & 7,056,461 & 7,058,781 & 7,090,320 & 0.44 & 124 \\
\texttt{s13\_t1351} & 6,856,205 & 6,863,829 & 6,915,240 & 0.74 & 120 & 6,247,774 & 6,249,757 & 6,271,560 & 0.35 & 112 \\
\texttt{s13\_t1372} & 8,350,945 & 8,361,823 & 8,459,580 & 1.16 & 146 & 7,452,535 & 7,454,846 & 7,474,860 & 0.27 & 135 \\
\texttt{s13\_t1532} & 9,050,466 & 9,056,799 & 9,115,680 & 0.65 & 154 & 8,128,086 & 8,130,272 & 8,160,240 & 0.37 & 146 \\
\texttt{s13\_t1591} & 8,716,345 & 8,720,590 & 8,757,720 & 0.42 & 153 & 8,246,490 & 8,248,519 & 8,269,680 & 0.26 & 152 \\
\texttt{s13\_t1598} & 8,590,185 & 8,596,695 & 8,664,360 & 0.78 & 149 & 7,670,143 & 7,672,742 & 7,706,700 & 0.44 & 139 \\
\texttt{s14\_t1500} & 8,835,895 & 8,845,046 & 8,880,420 & 0.40 & 154 & 8,222,074 & 8,224,192 & 8,245,860 & 0.26 & 150 \\
\texttt{s14\_t1511} & 9,167,303 & 9,173,645 & 9,257,400 & 0.90 & 157 & 7,817,050 & 7,818,507 & 7,838,760 & 0.26 & 144 \\
\texttt{s14\_t1629} & 10,146,016 & 10,151,858 & 10,212,780 & 0.60 & 177 & 8,333,186 & 8,335,086 & 8,358,000 & 0.27 & 149 \\
\texttt{s14\_t1755} & 10,283,388 & 10,289,863 & 10,356,000 & 0.64 & 176 & 9,408,259 & 9,409,947 & 9,432,780 & 0.24 & 170 \\
\texttt{s15\_t1713} & 10,034,985 & 10,041,287 & 10,130,340 & 0.88 & 172 & 8,823,104 & 8,825,561 & 8,869,920 & 0.50 & 156 \\
\texttt{s15\_t1745} & 10,366,864 & 10,372,707 & 10,458,720 & 0.82 & 177 & 9,112,295 & 9,114,199 & 9,140,520 & 0.29 & 164 \\
\texttt{s15\_t1749} & 9,308,113 & 9,312,168 & 9,372,480 & 0.64 & 162 & 8,622,494 & 8,624,236 & 8,653,620 & 0.34 & 158 \\
\texttt{s17\_t1955} & 11,758,675 & 11,766,551 & 11,854,920 & 0.75 & 204 & 10,117,389 & 10,119,161 & 10,152,660 & 0.33 & 186 \\
\bottomrule
\end{tabular}
\end{adjustbox}
\end{table}

%% file: tables/tableBestKnownLarge.tex
\begin{table}[htbp]
\centering
\caption{Best known bounds and solutions for large instances.}
\label{tab:best_large}
\begin{adjustbox}{max width=\textwidth}
\begin{tabular}{l rrr r r rrr r r}
\toprule
& \multicolumn{5}{c}{Single-depot} & \multicolumn{5}{c}{Multi-depot} \\
\cmidrule(lr){2-6} \cmidrule(lr){7-11}
Instance & Root LB & Best LB & Best UB & Gap (\%) & Duties & Root LB & Best LB & Best UB & Gap (\%) & Duties \\
\midrule
\texttt{s15\_t2029} & 11,134,303 & 11,140,680 & 11,257,260 & 1.04 & 190 & 9,917,645 & 9,919,298 & 9,960,660 & 0.42 & 177 \\
\texttt{s16\_t2052} & 11,938,880 & 11,944,278 & 12,045,540 & 0.84 & 203 & 10,468,828 & 10,470,472 & 10,514,220 & 0.42 & 187 \\
\texttt{s16\_t2136} & 12,715,984 & 12,719,602 & 12,788,520 & 0.54 & 217 & 11,405,783 & 11,406,933 & 11,446,920 & 0.35 & 209 \\
\texttt{s16\_t2219} & 13,318,246 & 13,324,371 & 13,437,960 & 0.85 & 227 & 10,702,291 & 10,703,286 & 10,741,800 & 0.36 & 195 \\
\texttt{s16\_t2252} & 12,396,038 & 12,400,856 & 12,535,440 & 1.07 & 214 & 11,064,280 & 11,065,751 & 11,118,060 & 0.47 & 200 \\
\texttt{s17\_t2090} & 12,436,106 & 12,441,737 & 12,525,420 & 0.67 & 212 & 10,852,090 & 10,853,424 & 10,891,680 & 0.35 & 197 \\
\texttt{s17\_t2166} & 14,037,157 & 14,044,780 & 14,128,740 & 0.59 & 238 & 11,642,240 & 11,643,649 & 11,681,700 & 0.33 & 214 \\
\texttt{s17\_t2183} & 14,290,908 & 14,297,875 & 14,451,960 & 1.07 & 244 & 12,328,601 & 12,331,002 & 12,381,840 & 0.41 & 219 \\
\texttt{s17\_t2214} & 12,387,240 & 12,391,876 & 12,441,960 & 0.40 & 216 & 10,657,805 & 10,659,956 & 10,720,080 & 0.56 & 187 \\
\texttt{s17\_t2242} & 14,742,098 & 14,748,939 & 14,952,300 & 1.36 & 256 & 12,606,446 & 12,608,013 & 12,650,520 & 0.34 & 229 \\
\texttt{s17\_t2366} & 13,478,416 & 13,483,797 & 13,620,840 & 1.01 & 234 & 11,891,063 & 11,892,482 & 11,931,600 & 0.33 & 218 \\
\texttt{s18\_t2395} & 15,425,380 & 15,429,645 & 15,542,940 & 0.73 & 263 & 13,764,435 & 13,765,215 & 13,799,940 & 0.25 & 249 \\
\texttt{s18\_t2550} & 14,939,689 & 14,944,106 & 15,041,760 & 0.65 & 257 & 13,480,955 & 13,481,910 & 13,541,220 & 0.44 & 243 \\
\texttt{s18\_t2595} & 13,723,693 & 13,728,804 & 13,837,860 & 0.79 & 234 & 12,149,343 & 12,150,418 & 12,199,080 & 0.40 & 223 \\
\texttt{s18\_t2641} & 15,474,385 & 15,478,417 & 15,627,720 & 0.96 & 265 & 13,458,876 & 13,459,831 & 13,528,980 & 0.51 & 237 \\
\texttt{s18\_t2648} & 13,801,297 & 13,806,594 & 13,921,500 & 0.83 & 234 & 12,094,401 & 12,095,027 & 12,150,300 & 0.45 & 218 \\
\texttt{s19\_t2315} & 14,810,475 & 14,817,684 & 14,958,960 & 0.94 & 253 & 12,585,061 & 12,586,784 & 12,633,420 & 0.37 & 225 \\
\texttt{s19\_t2484} & 15,219,679 & 15,225,464 & 15,369,840 & 0.94 & 261 & 12,461,767 & 12,463,222 & 12,537,360 & 0.59 & 226 \\
\texttt{s19\_t2661} & 18,565,600 & 18,571,086 & 18,741,660 & 0.91 & 317 & 14,820,822 & 14,821,928 & 14,882,040 & 0.40 & 269 \\
\texttt{s19\_t2702} & 14,331,231 & 14,335,649 & 14,423,520 & 0.61 & 247 & 12,840,737 & 12,840,737 & 12,911,940 & 0.55 & 233 \\
\texttt{s19\_t2711} & 15,915,564 & 15,919,393 & 15,988,620 & 0.43 & 270 & 14,476,229 & 14,477,151 & 14,533,020 & 0.38 & 258 \\
\texttt{s19\_t2744} & 15,100,876 & 15,104,590 & 15,214,140 & 0.72 & 259 & 13,331,065 & 13,331,810 & 13,382,700 & 0.38 & 243 \\
\texttt{s20\_t2668} & 15,799,233 & 15,802,963 & 15,861,120 & 0.37 & 274 & 13,605,974 & 13,606,853 & 13,665,660 & 0.43 & 253 \\
\texttt{s20\_t2921} & 18,007,663 & 18,011,628 & 18,133,260 & 0.67 & 310 & 15,992,132 & 15,993,143 & 16,076,100 & 0.52 & 291 \\
\texttt{s20\_t3016} & 16,181,075 & 16,184,454 & 16,287,780 & 0.63 & 277 & 14,703,742 & 14,704,066 & 14,756,220 & 0.35 & 265 \\
\bottomrule
\end{tabular}
\end{adjustbox}
\end{table}

%% file: figs/parameterTuning.tex
\begin{figure}[htbp]
\centering
\tikzset{
  legend matrix/.style={anchor=north, column sep=3pt, nodes={anchor=center}},
  legend mark/.style={draw=black, inner sep=1.5pt},
  legend label/.style={font=\small, inner sep=0pt, text height=1.5ex, text depth=0.25ex}
}
\begin{subfigure}[b]{0.48\textwidth}
\centering
\begin{tikzpicture}
\begin{axis}[
    width=\textwidth,
    xlabel={Average Number of Iterations},
    ylabel={Average Time (s)},
    tick label style={font=\tiny},
    label style={font=\small},
    mark size=1.5pt
]
\addplot[only marks, mark=*, mark options={fill=black, draw=black}] coordinates {
    (138.6000,65.6424)
    (122.5600,45.6257)
    (106.4800,42.4490)
    (123.1200,62.8613)
    (117.8800,45.6344)
    (105.8000,43.1964)
    (122.4400,67.3524)
    (115.6800,48.4862)
    (106.0000,44.4691)
};
\addplot[only marks, mark=square*, mark options={fill=gray!75, draw=black}] coordinates {
    (95.0800,49.3900)
    (100.2000,38.3101)
    (89.4000,36.4557)
    (84.2800,47.6951)
    (98.8000,37.6957)
    (86.1600,37.5388)
    (84.0000,53.0812)
    (92.5200,41.9284)
    (88.8000,39.3581)
};
\addplot[only marks, mark=diamond*, mark options={fill=white, draw=black}] coordinates {
    (86.8400,48.2826)
    (96.9600,36.7061)
    (85.4800,36.5675)
    (78.8800,48.3384)
    (91.4000,36.8062)
    (82.5600,37.2400)
    (76.8800,57.1424)
    (91.5200,41.1533)
    (88.6000,40.9912)
};
\end{axis}
\matrix[legend matrix] at ([yshift=-3pt]current axis.south |- current axis.outer south) {
  \node[legend mark, circle, fill=black] {}; & \node[legend label] {0.25}; &[0.5cm]
  \node[legend mark, rectangle, fill=gray!75] {}; & \node[legend label] {0.5}; &[0.5cm]
  \node[legend mark, diamond, fill=white] {}; & \node[legend label] {0.75}; \\
};
\end{tikzpicture}
\caption{Column selection threshold}
\end{subfigure}
\hfill
\begin{subfigure}[b]{0.48\textwidth}
\centering
\begin{tikzpicture}
\begin{axis}[
    width=\textwidth,
    xlabel={Average Number of Iterations},
    ylabel={Average Time (s)},
    tick label style={font=\tiny},
    label style={font=\small},
    mark size=1.5pt
]
\addplot[only marks, mark=*, mark options={fill=black, draw=black}] coordinates {
    (138.6000,65.6424)
    (122.5600,45.6257)
    (106.4800,42.4490)
    (95.0800,49.3900)
    (100.2000,38.3101)
    (89.4000,36.4557)
    (86.8400,48.2826)
    (96.9600,36.7061)
    (85.4800,36.5675)
};
\addplot[only marks, mark=square*, mark options={fill=gray!75, draw=black}] coordinates {
    (123.1200,62.8613)
    (117.8800,45.6344)
    (105.8000,43.1964)
    (84.2800,47.6951)
    (98.8000,37.6957)
    (86.1600,37.5388)
    (78.8800,48.3384)
    (91.4000,36.8062)
    (82.5600,37.2400)
};
\addplot[only marks, mark=diamond*, mark options={fill=white, draw=black}] coordinates {
    (122.4400,67.3524)
    (115.6800,48.4862)
    (106.0000,44.4691)
    (84.0000,53.0812)
    (92.5200,41.9284)
    (88.8000,39.3581)
    (76.8800,57.1424)
    (91.5200,41.1533)
    (88.6000,40.9912)
};
\end{axis}
\matrix[legend matrix] at ([yshift=-3pt]current axis.south |- current axis.outer south) {
  \node[legend mark, circle, fill=black] {}; & \node[legend label] {10}; &[0.5cm]
  \node[legend mark, rectangle, fill=gray!75] {}; & \node[legend label] {20}; &[0.5cm]
  \node[legend mark, diamond, fill=white] {}; & \node[legend label] {50}; \\
};
\end{tikzpicture}
\caption{Column removal threshold}
\end{subfigure}
\\[1em]
\begin{subfigure}[b]{0.48\textwidth}
\centering
\begin{tikzpicture}
\begin{axis}[
    width=\textwidth,
    xlabel={Average Number of Iterations},
    ylabel={Average Time (s)},
    tick label style={font=\tiny},
    label style={font=\small},
    mark size=1.5pt
]
\addplot[only marks, mark=*, mark options={fill=black, draw=black}] coordinates {
    (122.5600,45.6257)
    (117.8800,45.6344)
    (115.6800,48.4862)
    (100.2000,38.3101)
    (98.8000,37.6957)
    (92.5200,41.9284)
    (96.9600,36.7061)
    (91.4000,36.8062)
    (91.5200,41.1533)
};
\addplot[only marks, mark=square*, mark options={fill=gray!75, draw=black}] coordinates {
    (106.4800,42.4490)
    (105.8000,43.1964)
    (106.0000,44.4691)
    (89.4000,36.4557)
    (86.1600,37.5388)
    (88.8000,39.3581)
    (85.4800,36.5675)
    (82.5600,37.2400)
    (88.6000,40.9912)
};
\addplot[only marks, mark=diamond*, mark options={fill=white, draw=black}] coordinates {
    (138.6000,65.6424)
    (123.1200,62.8613)
    (122.4400,67.3524)
    (95.0800,49.3900)
    (84.2800,47.6951)
    (84.0000,53.0812)
    (86.8400,48.2826)
    (78.8800,48.3384)
    (76.8800,57.1424)
};
\end{axis}
\matrix[legend matrix] at ([yshift=-3pt]current axis.south |- current axis.outer south) {
  \node[legend mark, circle, fill=black] {}; & \node[legend label] {$\alpha=0.5$}; &[0.5cm]
  \node[legend mark, rectangle, fill=gray!75] {}; & \node[legend label] {Self-tuning $\alpha$}; &[0.5cm]
  \node[legend mark, diamond, fill=white] {}; & \node[legend label] {$\alpha=0$}; \\
};
\end{tikzpicture}
\caption{Dual smoothing setting}
\end{subfigure}
\caption{Effect of column selection, column management, and dual smoothing settings on the average number of pricing iterations and average root node solve time. Each marker represents one of the 27 parameter configurations, averaged over medium-sized instances.}
\label{fig:parameter_tuning}
\end{figure}